\documentclass[leqno,a4paper,oneside]{amsart}
\usepackage{amsmath,amssymb,amscd,mathrsfs}
\usepackage{graphicx,subfigure,caption,comment,multirow,color}
\usepackage{wrapfig}
\usepackage{mathptmx,courier}
\usepackage[scaled=1.0]{helvet}
\title{Strongly positive amphicheiral knots with doubly symmetric diagrams}

\author{Christoph Lamm}

\theoremstyle{plain}
  \newtheorem{theorem}{Theorem}[section]

\theoremstyle{definition}
  \newtheorem{definition}[theorem]{Definition}
  \newtheorem{question}[theorem]{Question}
	
  \newtheorem{remark}[theorem]{Remark}

\newcommand{\myS}{\mathbb{S}}

\begin{document} 

\begin{abstract}
We determine the prime strongly positive amphicheiral knots up to 16 crossings and show 
that a large fraction of them admit knot diagrams with a double symmetry (rotational symmetry 
for strongly positive amphicheirality and an additional mirror symmetry for the ribbon property). 
The remaining knots are presented as `almost doubly symmetric' diagrams, defined as diagrams
with double symmetry where exactly one symmetric pair of crossings is switched.
We also find that the rosette knot 14a19470, which has period 7, is the first prime strongly 
positive amphicheiral knot in the knot tables which is not slice.
\end{abstract}



\captionsetup{belowskip=10pt,aboveskip=10pt}

\def\scaling{0.8}

\definecolor{myboxcolour}{gray}{0.8}

\reversemarginpar


\maketitle

\section{Introduction} \label{sec:Introduction}

In the article on the tabulation of strongly invertible knots \cite{Lamm2022} we announced a study
on strongly positive amphicheiral knots. The current article contains the results of this study.
In particular, we list all prime strongly positive amphicheiral knots up to 16 crossings and present
diagrams with additional symmetries for them.

\begin{figure}[hbtp]
\centering
\includegraphics[scale=0.72]{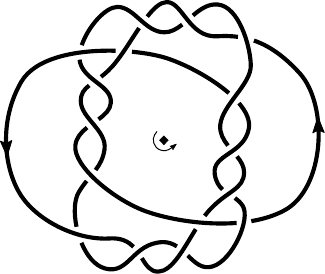}
\hspace{1.7cm}
\raisebox{-0.2cm}{\includegraphics[scale=0.8]{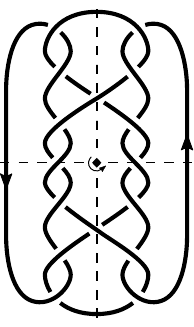}}
\caption{A rotationally symmetric diagram (for $14a_{19517}$) and a doubly symmetric diagram (for $12n_{706}$)}
\label{14a19517_12n706}
\end{figure}

A knot is called strongly positive amphicheiral if it has a diagram which is mapped to its mirror image 
by a rotation of $\pi$, preserving the orientation. An example is the diagram on the left in Figure \ref{14a19517_12n706}.

In many respects, strongly positive amphicheiral knots are similar to slice knots:
Their Alexander polynomial is always a square \cite{HartleyKawauchi} and this is similar to the $f(t)\cdot f(1/t)$ 
condition for the Alexander polynomial of slice knots (compare also with Long's result on algebraic sliceness \cite{Long}).
However, Flapan showed in 1986 that not all prime strongly positive amphicheiral knots are slice \cite{Flapan}.

\enlargethispage{1.0cm}

On the other hand, the family of strongly positive amphicheiral knots is quite smaller than that of slice knots.
For instance, 2-bridge knots cannot be strongly positive amphicheiral \cite{HartleyKawauchi}.

In the current study, we analyze diagrams of strongly positive amphicheiral knots which are also symmetric unions,
a property which implies that a knot is slice:

In \cite{Lamm2021} we found that all prime strongly positive amphicheiral knots with crossing number up to 12 possess diagrams 
with additional symmetries. This type of diagram with additional symmetries is simultaneously a symmetric union (with respect to a 
vertical axis) and has strongly invertible partial knots (with respect to a horizontal axis of rotation in a transvergent diagram). 

This combination of two simultaneous symmetries results again in a rotationally symmetric (strongly positive amphicheiral) 
diagram with respect to an axis perpendicular to the diagram plane, see Section \ref{sec:diagrams} for more details. 
We call these diagrams `doubly symmetric'.

\begin{figure}[hbtp]
\centering
\includegraphics[scale=0.75]{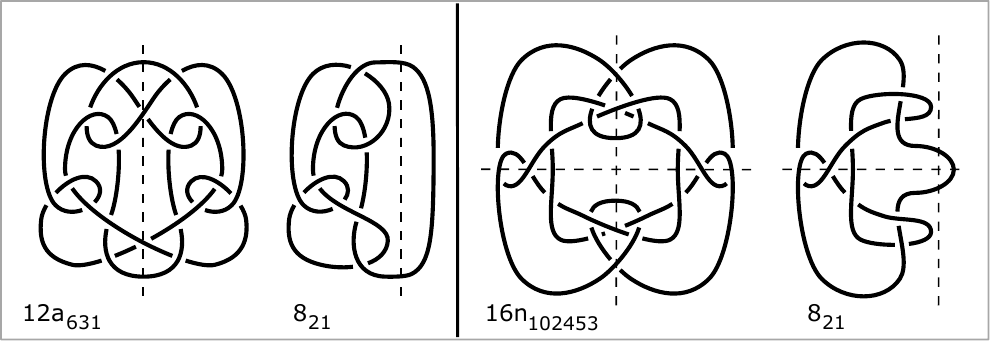}
\caption{The partial knot in a general symmetric union diagram (left) and in a doubly symmetric diagram (right)}
\label{overview}
\end{figure}

In Figure \ref{overview} a general symmetric union with partial knot $8_{21}$ is contrasted with a symmetric union in a doubly 
symmetric diagram (with the same partial knot). The diagram on the right in Figure \ref{14a19517_12n706} is another example of
a doubly symmetric diagram.

Knots with symmetric union diagrams are ribbon knots. Therefore, knots with doubly symmetric diagrams combine 
properties of strongly positive amphicheiral knots and ribbon knots. Our initial question was:

\begin{question}
Which prime strongly positive amphicheiral knots up to 16 crossings have doubly symmetric diagrams?
\end{question}

Using SnapPy \cite{SnapPy}, it was possible to determine all prime strongly positive amphicheiral knots up to 16 crossings.
There are 131 of them and they are listed in Tables \ref{tab:knotList} and \ref{tab:knotList2}.

Although we found that from these at least 106 have doubly symmetric diagrams, it is clear that not all prime strongly positive 
amphicheiral knots will have this property: As mentioned, not all prime strongly positive amphicheiral knots are ribbon.

Before continuing the summary of results, we dive into a non-mathematical digression.

\subsection{Doubly symmetric diagrams, knot gardens and ornaments}
I encountered the term `knot garden' in Andrew Ranicki's slides \cite{Ranicki}, where on page 6 a Tudor knot garden is mentioned.
The illustration on this slide reminded me of the symmetry of doubly symmetric knot diagrams and therefore I include a
digression on knot gardens and symmetric ornaments.

Knot gardens appeared in Europe around 1500. Richard M. Bacon (freelance writer 
and horticultural consultant) wrote in 1985 in the New York Times about them:

\enlargethispage{1.0cm}

\begin{quote}
The idea of weaving together lines of plants to establish this pattern came from monastic enclosures in medieval Europe. 
Herbs were most certainly used because they were cultivated to assure good health.
The device was taken over later by lay gardeners and made even more decorative. 
For a time they were at the height of fashion and became the focal point for classical and formal gardens in this country by 1800. 
Then they were expanded and superceded by the parterre which introduced internal paths to the pattern and a broader selection of plant material. 
\end{quote}

\begin{figure}[hbtp]
\centering
\includegraphics[scale=0.3]{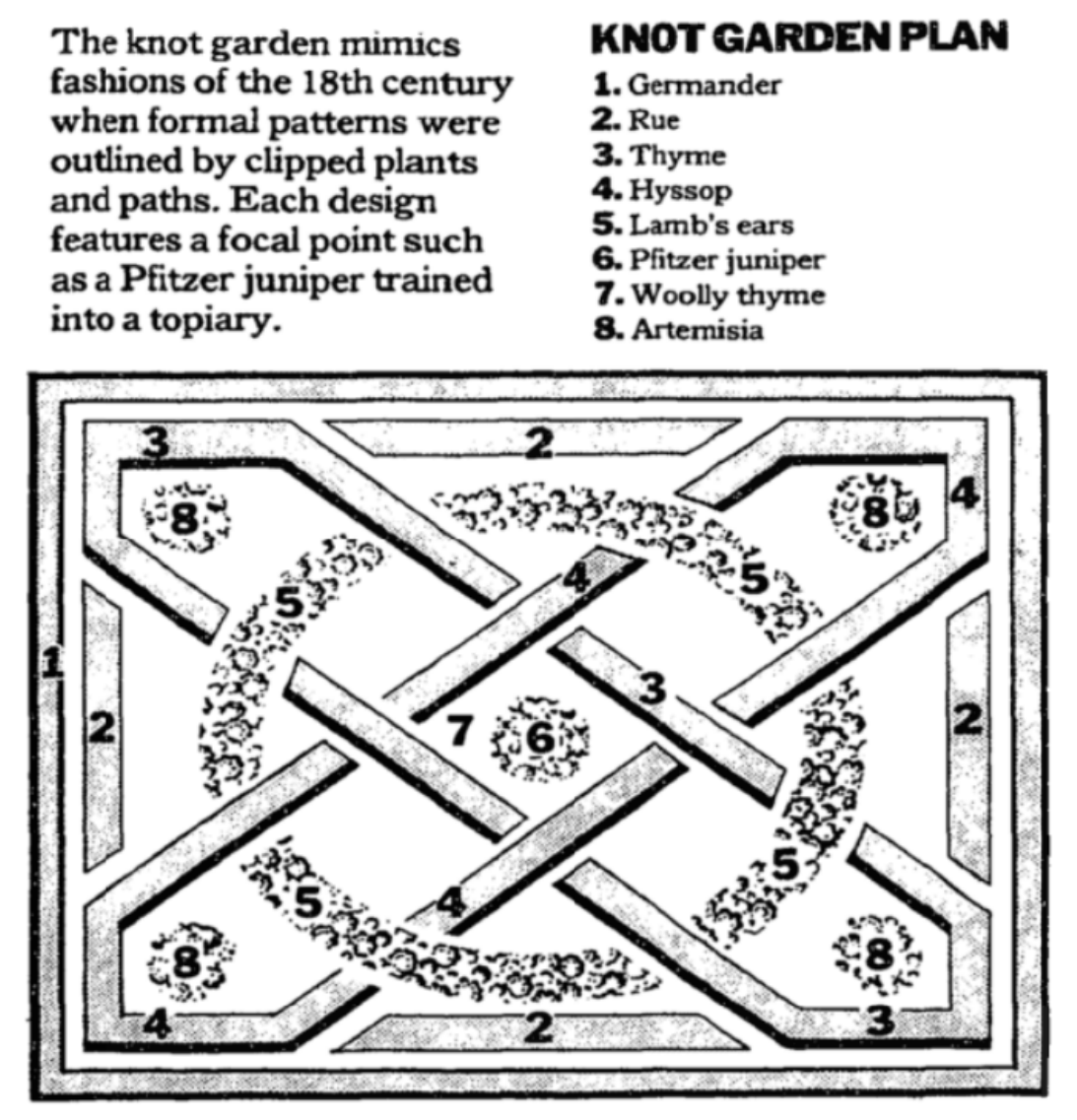}
\caption{A knot garden illustration, New York Times, May 19th, 1985}
\label{NYT_garden_plan}
\end{figure}

The NYT article contains the illustration in Figure \ref{NYT_garden_plan}. The garden contains four closed components and several additional elements.
A knot garden in a strict mathematical sense would contain one component only and no other decorative additions.

A related term, occuring in the last sentence, is the `parterre'. Parterres are similar if seen from above 
but they miss the three-dimensional weaving of knot gardens. As a reminiscence to horticulture, we will use 
the term parterre in this article for knot shadows, i.e. knot diagrams without over- and under-crossings.
The concrete meaning of these terms for historic gardens in England seems to be quite unclear; see the appendix
on knot gardens in \cite{Francis}.

As an example for a knot garden in literature we found the book `Hypnerotomachia Poliphili' (The Dream of Poliphilus), which appeared
in 1499 and is attributed to the Italian monk Francesco Colonna. An analysis of the knot garden theme in this work can be found in \cite{Segre}.
There is even an opera with the title `The knot garden' by the English composer Michael Tippett  (first performed in 1970).


We mention that the required symmetry for the parterre (gardening and doubly symmetric diagrams) is the dihedral symmetry of a rectangle ($D_2$).
The additional symmetry of a square is not needed, although in many cases a garden compartment will be a square and might have a $D_4$ symmetry.

\enlargethispage{0.5cm}

\begin{figure}[hbtp]
\centering
\includegraphics[scale=0.23]{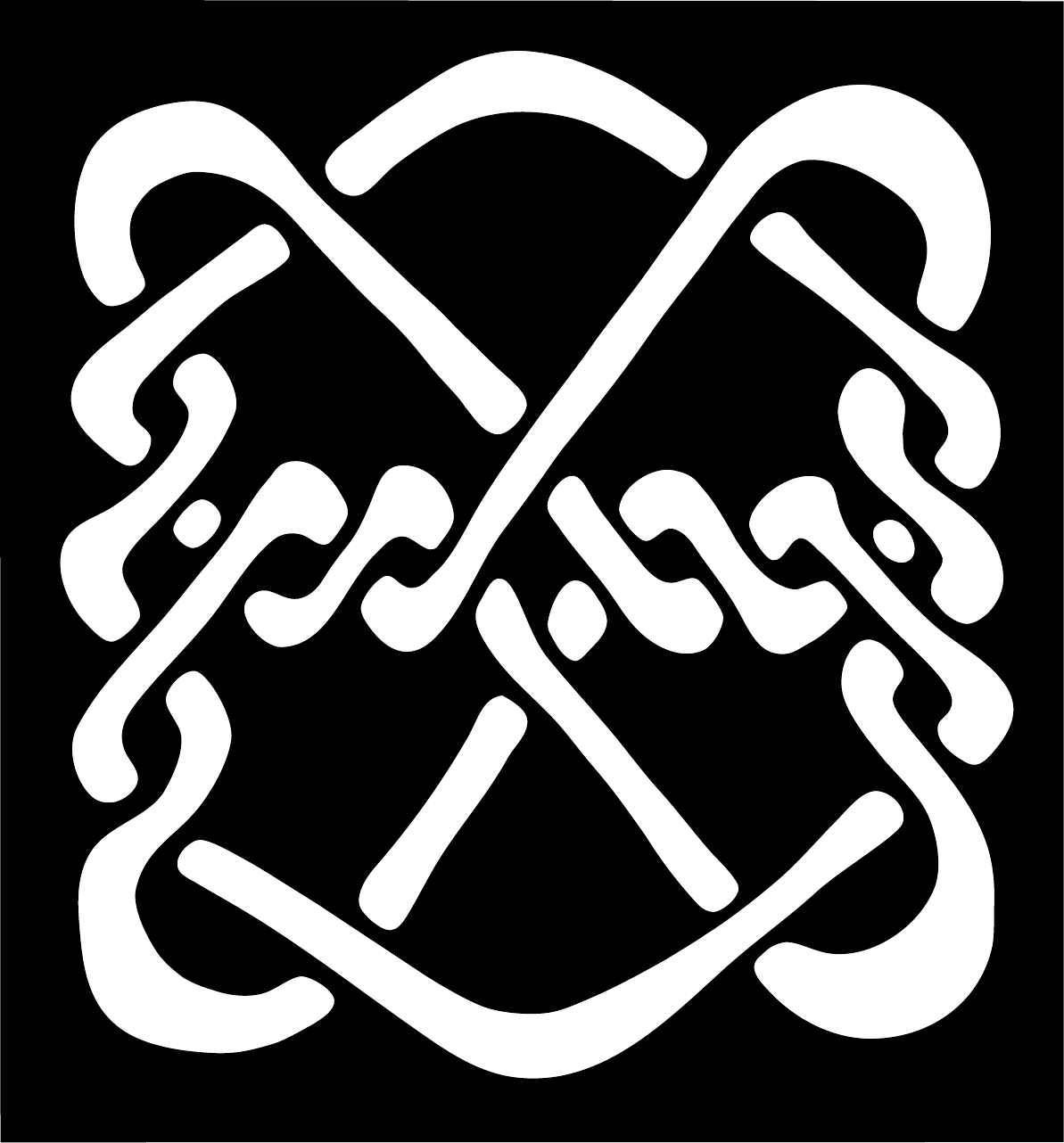}
\quad
\includegraphics[scale=0.21]{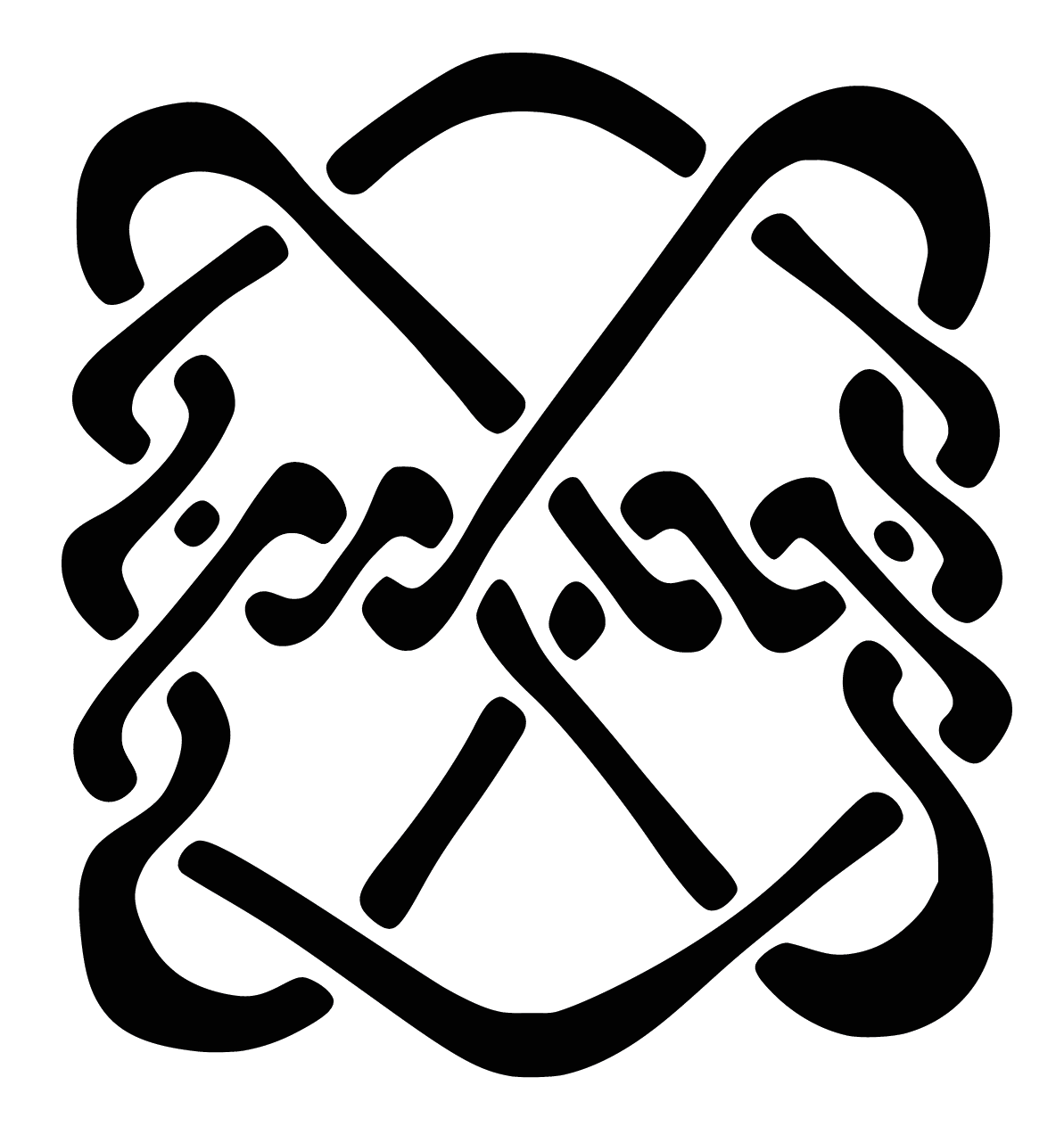}
\caption{A lino print of $16n_{872167}$ (left), and an inverted version (right), showing the diagram $t_4(0, 0, 3 \mid 1)$ in template notation}
\label{linoprint_16n872167}
\end{figure}

\newpage 

Symmetric designs also occur in celtic and islamic ornaments. As an example which is close to a celtic 
knot design, we did a linoleum print of a doubly symmetric diagram, see Figure \ref{linoprint_16n872167}. 
Note, though, that in celtic knot design we almost never find non-alternating diagrams. 
The use of non-alternating diagrams cannot be avoided when knot diagrams are constructed as symmetric unions 
because for the arcs which traverse the vertical axis two under- or two over-crossings are necessarily adjacent.

In Versailles we find the following parterre with two symmetry axes ($D_2$-symmetry):

\begin{figure}[hbtp]
\centering
\includegraphics[scale=0.6]{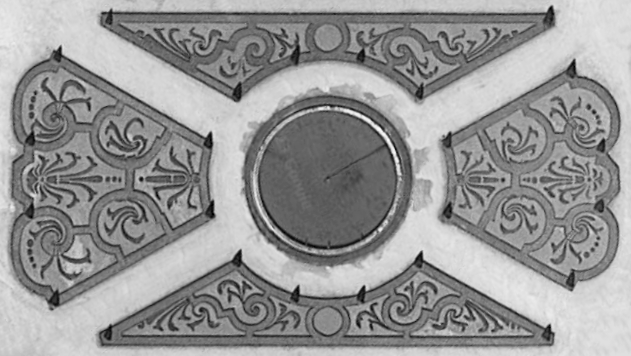}
\caption{Part of the `Parterre du Midi', Versailles}
\label{Parterre_du_midi}
\end{figure}

\subsection{Doubly symmetric diagrams with additional symmetries}
Before we give a further outline of the article, we present the following doubly symmetric knot diagrams with $D_4$-symmetry of their parterres.
The second diagram can be rotated by 45° (counter-clockwise) into another doubly symmetric diagram and therefore has two sets of axes simultaneously.
This is not the case, however, for the other two diagrams shown here. (As an exercise we invite the reader to check these statements.)

\begin{figure}[hbtp]
\centering
\includegraphics[scale=0.8]{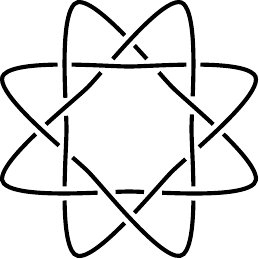}\hspace{0.4cm}
\includegraphics[scale=0.8]{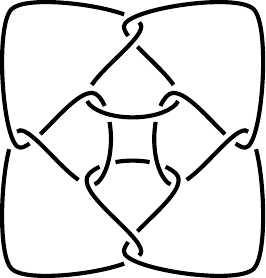}\hspace{0.4cm}
\includegraphics[scale=0.8]{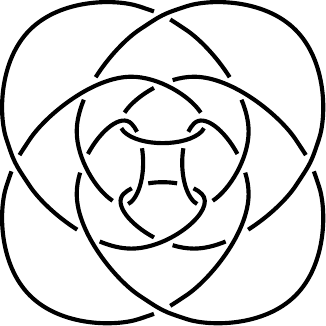}
\caption{Three doubly symmetric knot diagrams for $14a_{19472}$, $16n_{991381}$ and $16n_{268004}$ with additional diagonal symmetries.
These diagrams are $t_4(0, 1, 0 \mid 1)$, $t_1(2, 0 \mid 2)$ and $t_9(1, 1, 0 \mid -1, 1)$.}
\label{additional_sym}
\end{figure}

\newpage

\subsection{Plan of the article}
We continue the overview of the results and the article's structure. 
In the next section we explain how, with the help of SnapPy, we determined the 131 prime strongly positive amphicheiral knots up to 16 crossings.
We found experimentally that from these at least 106 have doubly symmetric diagrams and therefore are ribbon knots.
For the remaining 25 knots we found `almost doubly symmetric' diagrams. These are diagrams with double symmetry where exactly one 
symmetric pair of crossings is switched. They are no longer symmetric unions, but still represent strongly positive amphicheiral knots.

\begin{question}
Does every strongly positive amphicheiral knot have a doubly symmetric or almost doubly symmetric diagram?
\label{question_doubly_almost}
\end{question}

The main open problem is to find a knot invariant which enables us to show that there is a strongly positive 
amphicheiral knot which is ribbon, or even a symmetric union, but does not have a doubly symmetric diagram.

\section{Strongly positive amphicheiral knots up to 16 crossings}
Surprisingly, prime strongly positive amphicheiral knots up to 16 crossings have not yet been tabulated.
Before we describe the SnapPy script which does the tabulation, we explain in an example how SnapPy is used for this purpose.

Until recently it was unclear whether the knot $12a_{435}$ is strongly positive amphicheiral (it is an amphicheiral knot 
and its Alexander polynomial is a square). For the analysis of this knot we proceed as follows with SnapPy:

\verb/S=Manifold('12a435').symmetry_group()/: SnapPy has stored all hyperbolic knots up to 16 crossings and the hyperbolic 
structures of their complements. It either computes the symmetry group (of isometries) on the spot or it is also already stored.
The variable $S$ now contains this symmetry group.

\verb/S.isometries()/: Gives the list of elements of the symmetry group. In this case there are 8 isometries and the output looks as follows (we omit the information `Extends to link' and the component mapping \verb/0 -> 0/ for each of the isometries):

\small
\begin{verbatim}
 [1 0]  [-1 0]  [1 0]  [-1 0]  [1  0]  [-1  0]  [1  0]  [-1  0] 
 [0 1]  [ 0 1]  [0 1]  [ 0 1]  [0 -1]  [ 0 -1]  [0 -1]  [ 0 -1]
\end{verbatim}
\normalsize

Explanation: The matrix refers to the meridian and longitude mapping. The first of the eight isometries is the identity.
The isometries are in this case numbered from 0 through 7. Their effect on the orientation of $\myS^3$ and on the knot is as follows:

\small
\medskip

\noindent
Isometries 0 and 2 are orientation-preversing isometries preserving the knot's orientation.  \newline
Isometries 1 and 3 are orientation-reversing isometries preserving the knot's orientation.   \newline
Isometries 4 and 6 are orientation-reversing isometries reversing the knot's orientation.    \newline
Isometries 5 and 7 are orientation-preversing isometries reversing the knot's orientation.
\normalsize

\medskip
The knot $12a_{435}$ is positive amphicheiral (isometries 1 and 3) and negative amphicheiral (isometries 4 and 6).
If one of the isometries 1 or 3 is an involution, the knot is strongly positive amphicheiral.

Multiplication of isometry 1 with itself is done with \verb/S.multiply_elements(1,1)/: The result is 2, which is not the identity. 
Therefore isometry 1 is not an involution.
We do the same for isometry 3, and the result is again 2. Hence the knot $12a_{435}$ is not strongly positive amphicheiral.

The following script does this for all amphicheiral knots in a specified range. It is called for alternating and non-alternating 
knots with crossing number 10, 12, 14 and 16 (line 3 is modified for that and the range in line 2 has to be adapted as well):

\newpage

\small
\begin{verbatim}
 1 def print_amphicheiral():
 2     for n in range(1,888):
 3         K_id = '12n'+str(n)
 4         K = Manifold(K_id)
 5         try:
 6             sym = K.symmetry_group()
 7         except:
 8             print(K_id, " is not hyperbolic.")
 9         if sym.is_amphicheiral():
10             iso = sym.isometries()
11             pa  = false  # positive amphicheiral flag
12             spa = false  # strongly positive amphicheiral flag
13             count = 0
14             for i in iso:
15                 A = i.cusp_maps()[0]
16                 if A[0,0] == -1 and A[1,1] == 1:
17                     pa = true
18                     if sym.multiply_elements(count,count) == 0:
19                         spa = true
20                 count += 1
21             if spa:
22                 print('knot ', K_id, ' is strongly positive amphicheiral')
23             elif pa:
24                 print('knot ', K_id, ' is positive amphicheiral')
25             else:
26                 print('knot ', K_id, ' is amphicheiral')
\end{verbatim}

\normalsize
\medskip
To facilitate the understanding, we give more comments on individual code lines: For non-hyperbolic knots an exception is caught in line 7
(the symmetry group cannot be calculated by SnapPy). Line 9 filters for amphicheiral knots; the following lines are analogous to
our example and two boolean flags are defined. Lines 14 through 20 contain the iteration through the isometries. The matrix elements
are tested in line 16 and the involution property in line 18.

The result consists of 131 prime knots; they are listed in Tables \ref{tab:knotList} and \ref{tab:knotList2}.
This list is complete because prime amphicheiral knots up to 16 crossings are always hyperbolic:
Torus knots are always chiral and satellite knots up to 16 crossings are satellites of the trefoil 
and therefore chiral (see \cite{HosteThistlethwaiteWeeks}, p. 43).

\section{Templates for doubly symmetric diagrams} \label{sec:diagrams}

In this section, we describe the templates for tabulating doubly symmetric diagrams.
In the article \cite{Lamm2021} we already used templates for symmetric unions and \textsl{complete} diagrams were shown.
In \cite{Lamm2022} the table for strongly invertible knots used templates with only \textsl{half} of a diagram.
Now, because of the double symmetry, we may even reduce this to one quarter of a diagram!

\begin{figure}[hbtp]
\centering
\includegraphics[scale=0.65]{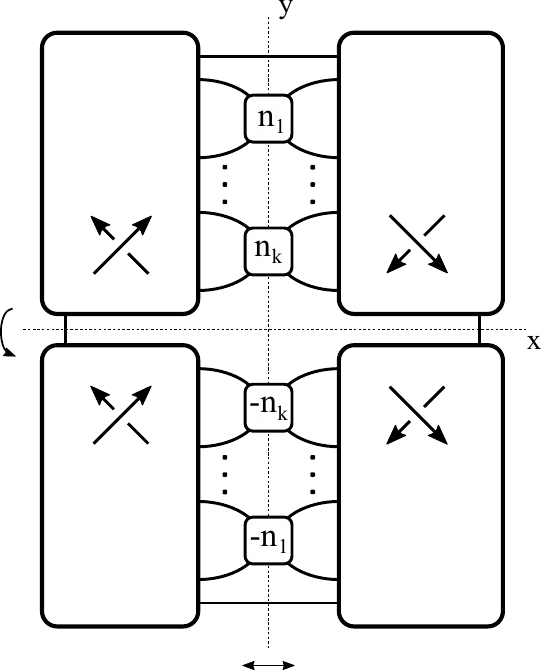}
\caption{The definition of doubly symmetric diagrams. The crossings on the axis $x$ are omitted in this illustration.
They are mirror symmetric with respect to the axis $y$ because doubly symmetric diagrams are symmetric unions.}
\label{doubly_symmetric_def}
\end{figure}

For doubly symmetric diagrams, we combine the properties of symmetric unions and strongly invertible (transvergent) diagrams.
These are defined as follows:
Symmetric unions are the connected sum of a diagram and its mirror image with reversed orientation (defining the \textsl{partial knot} 
of the symmetric union, which is unique up to mirror image), with additional crossings inserted on the symmetry axis. 
Transvergent diagrams for strongly invertible knots are rotationally symmetric with respect to an axis in the diagram plane, see \cite{Lamm2022} for details. 
These two constructions are therefore very similar and their parterre diagrams are indistinguishable. Care is needed for the choice of crossings, so that
the resulting diagram with two symmetries is indeed strongly positive amphicheiral:

\begin{definition}
A knot diagram $D$ is called \textit{doubly symmetric} if it possesses the two following symmetries with respect to perpendicular 
axes $x$ and $y$ in the diagram plane
\begin{itemize}
	\item[a)] $D$ is a symmetric union with axis $y$,
	\item[b)] the partial knot diagrams on the left and on the right sides are (mirror symmetric) transvergent diagrams with axis $x$,
\end{itemize}
and the inserted crossings on the axis $y$ are mirror symmetric with respect to $x$ as shown in Figure \ref{doubly_symmetric_def}.
A knot is called doubly symmetric if it has a doubly symmetric diagram.
\end{definition}

It is quite obvious that a knot with a doubly symmetric diagram is strongly positive amphicheiral:
The combination of rotation around the $x$-axis and reflection at the $y$-$z$-plane results in a point reflection
at the origin and the crossings on the two axes are chosen in a way that they are also point symmetric to the origin.
This is equivalent to rotating the diagram by $\pi$ into its mirror image and since the origin is not part of the diagram
we have a strongly \textit{positive} amphicheiral knot.

The only subtlety is the choice of the crossings on the $y$-axis. See the following remark.

\begin{remark}
We observe that the twist numbers for the crossings on the axis $y$ are $n_1, \ldots, n_k$ and $-n_k, \ldots, -n_1$.
A rotationally symmetric diagram with respect to $x$ would have twist number $n_1, \ldots, n_k$ and $n_k, \ldots, n_1$.
In this case, the total diagram would be a (strongly invertible) transvergent diagram, but it would lose the rotational symmetry 
with respect to the central axis perpendicular to the diagram plane which guarantees the strongly positive amphicheirality.
Therefore, although the partial knots are strongly invertible in this construction, doubly symmetric knots can be invertible or non-invertible.
\end{remark}

\begin{figure}[hbtp]
\centering
\includegraphics[scale=0.8]{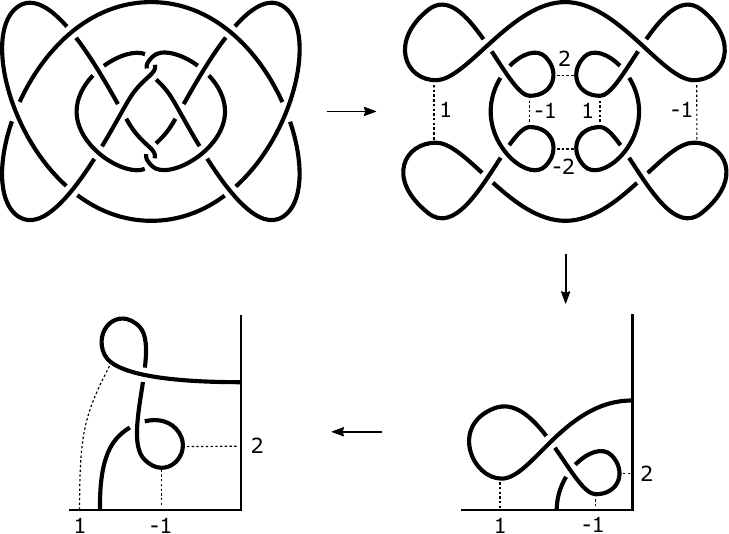}
\caption{Converting a doubly symmetric knot diagram into a template. The example shows a diagram of $16n_{645918}$, transformed into template $t_2 (1, -1 \mid 2)$.}
\label{16n645918_diagram_to_template}
\end{figure}

Our tabulating approach uses templates, similar to the way we did in \cite{Lamm2021} and \cite{Lamm2022}.
Each template contains the quarter diagram in the second quadrant and integer markings on the two axes for the twists which are inserted there. 
Figure \ref{16n645918_diagram_to_template} shows how a template diagram is constructed from  a doubly symmetric diagram.

\begin{figure}[hbtp]
\centering
\includegraphics[scale=1.0]{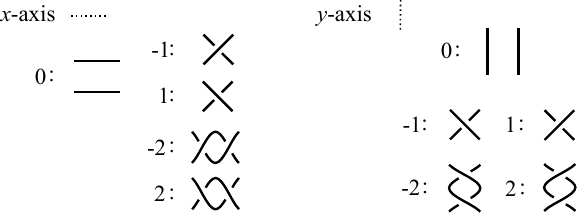}
\caption{The conventions for the integer markings for twists on the two axes.}
\label{twist_convention}
\end{figure}

\begin{figure}[hbtp]
\centering
\includegraphics[scale=0.8]{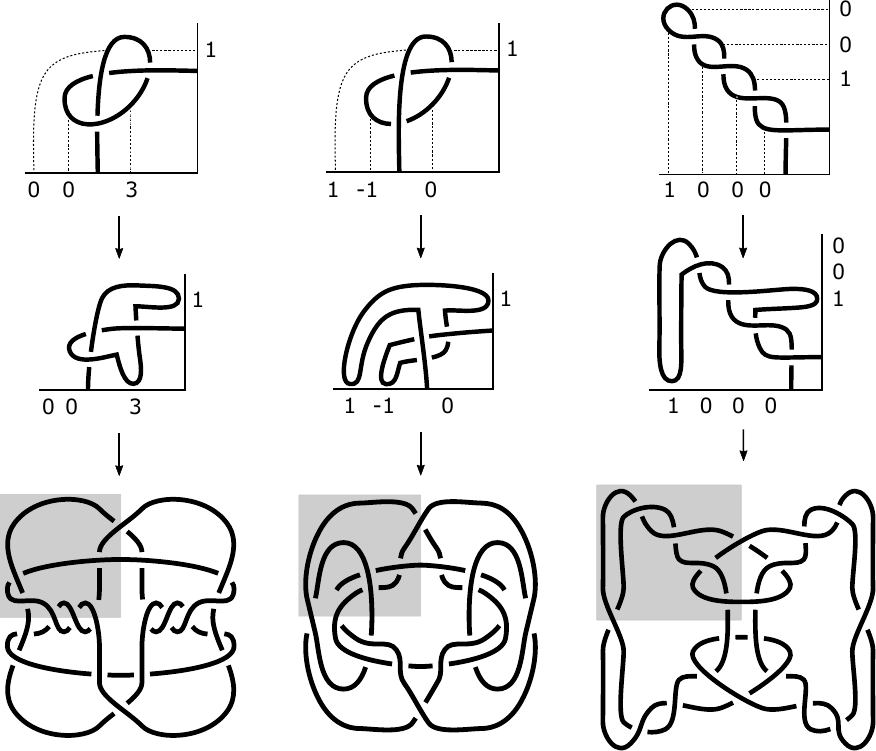}
\caption{Converting a template into a doubly symmetric knot diagram, illustrated by three doubly symmetric diagrams for the knot $16n_{872167}$:
  $t_4(0, 0, 3 \mid 1)$, $t_6(1, -1, 0 \mid 1)$ and $t_{17}(1, 0, 0, 0 \mid 0, 0, 1)$.}
\label{16n872167_3_constructions}
\end{figure}

In the other direction, we proceed as in Figure \ref{16n872167_3_constructions}: 
If a template diagram is given, we reconstruct the other quadrants and insert the twists on the axes given by the integer markings.
The conventions for these markings are defined in Figure \ref{twist_convention}.

The list of strongly positive amphicheiral knots up to 16 crossings in the appendix was compiled by starting with the simplest
templates made out of one or two crossings and by testing twist variations with Knotscape. We then added templates with three 
crossings, and for some knots templates with four crossings were necessary. In the appendix, we also list variations with different
partial knots because we consider doubly symmetric knot diagrams with different partial knots as interesting enough to be documented 
in detail. For instance, Figure \ref{16n872167_3_constructions} shows three different diagrams with partial knots $7_6$,
$8_{20}$ and $9_1$ for $16n_{872167}$. If several diagrams were found for a pair (consisting of a doubly symmetric knot and the partial knot),
we chose a diagrams with smallest crossing number for Table \ref{tab:knotList} and for the appendix.

\newpage

\section{Almost doubly symmetric diagrams}
\begin{figure}[hbtp]
\centering
\includegraphics[scale=0.8]{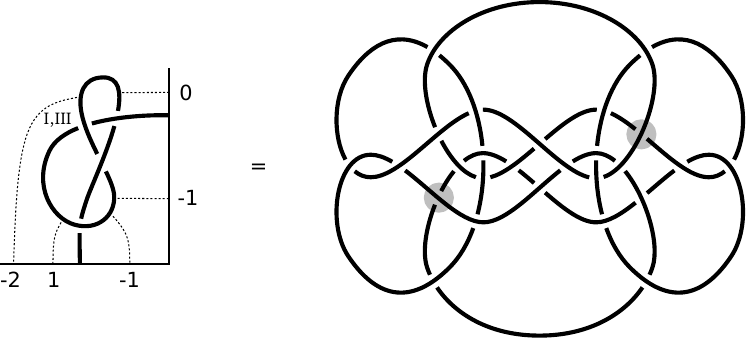}
\caption{An almost doubly symmetric diagram for $16a_{288139}$. The example shows the diagram constructed 
from template $t_{24}(-2, 1, -1 \mid 0, -1)$. The two switched rotationally symmetrical crossings are marked in gray.}
\label{almost_doubly_16a288139}
\end{figure}

The total number of prime strongly positive amphicheiral knots up to 16 crossings is 131 and from these we found 106 with doubly symmetric diagrams.
Wishing to show diagrams with a visible strongly positive amphicheiral symmetry also for the remaining knots,
we can choose between a) rotationally symmetric diagrams in braid form (see the diagram in the left of Figure \ref{14a19517_12n706}), or
b) rotationally symmetric diagrams built from more general tangles. We found by incidence that a small deviation from 
a doubly symmetric diagram yielded additional knots, and this led to the search for the remaining 25 knots as `almost
doubly symmetric' diagrams. Hence, we chose b) but with additional restrictions on the tangle form. Almost doubly symmetric diagrams
are doubly symmetric diagrams where two rotationally symmetrical crossings are modified, see Figure \ref{almost_doubly_16a288139}.
A knot is called almost doubly symmetric if it has an almost doubly symmetric diagram.

The template notation for almost doubly symmetric diagrams consists of an additional marking for one crossing: The quadrants in 
which the two symmetrical crossings are switched are indicated in the templates as {\small \verb/I,III/}, or {\small \verb/II,IV/}.
Currently it is not known whether this relatively small deviation from a double symmetry is enough to generate all
strongly positive amphicheiral knots, see Question \ref{question_doubly_almost}.

We remark, that due to the vertical symmetry it is possible to use only one of the markings, and for our templates we use {\small \verb/I,III/}.
If several rotationally symmetric crossing pairs are modified it can be necessary to apply both markings in one template description.

\newpage

\section{Summary data for prime amphicheiral knots up to 16 crossings}
Although the topic of this article is strongly positive amphicheiral knots, we would like to
give some data for amphicheiral knots in general. The total number of prime knots up to 16 crossings is 1,701,936 
(see the title of \cite{HosteThistlethwaiteWeeks}) if for each knot we count inverse and mirrored variants only once.
Amphicheiral knots contribute only with roughly $1/900 = 0.1111\%$ to this number, as there are 1892 of them (not including
the trivial knot). The distribution of amphicheiral knots into sub-classes is also interesting. 

The following Table \ref{tab:summary_data} splits the two symmetry classes of positive and fully amphicheiral knots into two parts, 
according to having the `strong' property (indicated by `,s'), and summarizes the numbers of amphicheiral knots.
The 131 prime strongly positive amphicheiral knots are contained in the columns `$+$,s' and `$a$,s', and there are 59 knots of type `$+$,s'
and 72 of type `$a$,s'. 

\begin{table}[htbp]
\begin{tabular}{|r|r|rr|rr|r|}
\hline
 c  & $-$ & $+$ & $+$,s & $a$ & $a$,s & Sum \\
\hline
 4  &      &    &    &  1 &     &   1 \\
 6  &      &    &    &  1 &     &   1 \\
 8  &   1  &    &    &  4 &     &   5 \\
10  &   6  &    &    &  5 &   2 &  13 \\
12  &  40  &    &  1 & 13 &   4 &  58 \\
14  & 227  &  2 &  4 & 25 &  16 & 274 \\
15  &   1  &    &    &    &     &   1 \\
16  & 1361 & 11 & 54 & 63 &  50 & 1539 \\
\hline
Sum & 1636 & 13 & 59 & 112 & 72 & 1892 \\
\hline 
\end{tabular}
\caption{The numbers of prime amphicheiral knots for each minimal crossing number $c$. 
The notation for the symmetry types is based on \cite{HosteThistlethwaiteWeeks}: 
`$-$' = negative amphicheiral, `$+$' = positive amphicheiral, `$a$' = fully amphicheiral, with `,s' added for strong symmetry.}
\label{tab:summary_data}
\end{table}


We give some comments on the sub-classes of amphicheiral knots:
By far the most frequent prime amphicheiral knots in this table are the negative amphicheiral knots. 
The only amphicheiral knot with 15 crossings (15n139717) is also of this type. The rarest knot family
consists of positive amphicheiral knots which are not strongly positive amphicheiral (13 knots in column `$+$').

We explicitly list the knots which occur 4 times or less (i.e. with table entries $\le 4$): 
The two amphicheiral knots with 4 and 6 crossings are $4_1$ and $6_3$. 
The knot with $c=8$ in column `$-$' is $8_{17}$. This is the first non-invertible knot in the knot tables.
The four fully amphicheiral knots with 8 crossings are $8_3$, $8_9$, $8_{12}$ and $8_{18}$.
The two strongly positive amphicheiral knots with 10 crossings are $10_{99}$ and $10_{123}$ (both invertible).
For $c=12$ we have the strongly positive amphicheiral knots 12a427 (non-invertible) and 12a1019, 12a1105, 12a1202, 12n706 (invertible).
For $c=14$ we have in column `$+$' the knots 14a10435 and 14n14148 (non-invertible) and in column `$+$,s' the 
knots 14a8662, 14a16309, 14a18676 and 14a18680 (also non-invertible).

\subsection{Double symmetry and invertibility}
We already mentioned that doubly symmetric knots may be invertible or non-invertible. The following
Table \ref{tab:doublesymmetry_invertibility} shows that the invertibility property is similarly distributed
for doubly and almost doubly symmetric knots. Note, that some of the knots which are recorded here as
almost doubly symmetric might be doubly symmetric.

\begin{table}[htbp]
\begin{tabular}{|l|r|r|r|}
\hline
  & $+$,s & $a$,s & Sum \\
\hline
doubly symmetric        & 45 & 61 & 106 \\
almost doubly symmetric & 14 & 11 &  25 \\
\hline
Sum                     & 59 & 72 & 131 \\
\hline 
\end{tabular}
\caption{The distribution of invertible (`$a$,s') and non-invertible (`$+$,s') knots among doubly and almost doubly symmetric knots.}
\label{tab:doublesymmetry_invertibility}
\end{table}

\subsection{Symmetry groups}
In Table \ref{tab:symmetry_groups} we give an overview of the symmetry groups of amphicheiral knots up to 16 crossings. 
This is based on Table A2A in \cite{HosteThistlethwaiteWeeks} and the result of our SnapPy script.
$C_n$ denotes the cyclic group with $n$ elements and $D_n$ the dihedral group with $2n$ elements (and therefore $C_2 = D_1$).

\begin{table}[htbp]
\begin{tabular}{|l|r|r|r|}
\hline
type  & most frequent case & additional cases \\
\hline
$-$   & $1577 \times D_1$ & $59 \times D_2$ \\
\hline
$+$   &   $13 \times C_4$ & \\
\hline
$+$,s &   $59 \times C_2$ & \\
\hline
$a$   &  $106 \times D_4$ & $4 \times D_8\;\,$    \\
      &                   & $1 \times D_{16}$ \\
\hline
$a$,s &   $66 \times D_2$ & $4 \times D_6\;\,$    \\
      &                   & $1 \times D_{10}$ \\
	    &                   & $1 \times D_{14}$ \\
\hline 
\end{tabular}
\caption{The symmetry groups occuring for each symmetry type}
\label{tab:symmetry_groups}
\end{table}

The fully amphicheiral knots with dihedral symmetry $D_{2n}$, $2n \ge 6$ are:
\begin{itemize}
\item[-] $D_6$: \: 12a1019, 12a1202, 12n706, 16n524438
\item[-] $D_8$: \; $8_{18}$, 16a357361, 16a375396, 16n995250
\item[-] $D_{10}$: $10_{123}$
\item[-] $D_{14}$: 14a19470
\item[-] $D_{16}$: 16a379778
\end{itemize}

\noindent
The last three of these knots are $n$-periodic because they are braid closures of the form $(\sigma_1 \sigma_2^{-1})^n$.
The knots with $c=12$ with symmetry group $D_6$ have period 3 \cite{JabukaNaik}, and $8_{18}$ has period 4. 

We look at the more detailed symmetry types according to \cite{BoyleRouseWilliams}.
For amphicheiral knots the following three types are possible: SNAP, RRef and SNASI.
Negative amphicheiral knots are of type SNAP. They have symmetry groups $D_n$ and period $n$.
Positive amphicheiral knots are of type RRef. They have symmetry groups $C_{2n}$ and period $n$.
Fully amphicheiral knots are of type SNASI, with symmetry groups $D_{2n}$ and period $n$.
In particular, 16n524438 also has period 3, and 16a357361, 16a375396, 16n995250 have period 4.
(Note, that hyperbolic amphicheiral knots are never freely periodic.)

\subsection{Ribbon property}
The knots in Table \ref{tab:knotList} are doubly symmetric and therefore ribbon. We checked the ribbon property 
of the 25 knots in Table \ref{tab:knotList2} with information taken from \cite{DunfieldGong} and \cite{OwensSwenton}. 
All knots but the following four are shown there to be ribbon knots: 14a19470, 16a288139, 16a289378, 16a309401.
The ribbon property of the last three of them is currently unknown, while 14a19470 is not slice \cite{Sartori}.
Therefore 14a19470 is the first knot in the knot tables which is strongly positive amphicheiral and not slice (compare \cite{Flapan}).

\begin{figure}[hbtp]
\centering
\includegraphics[scale=0.8]{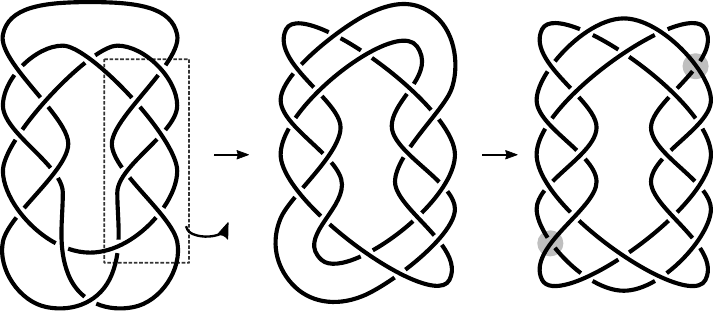}
\caption{A transformation of the knot $14a_{19470}$ from 7-periodic braid form into an almost doubly symmetric diagram, given by template $t_{24}(0, 0, -1 \mid 0, -1)$}
\label{almost_doubly_14a19470}
\end{figure}

In Figure \ref{almost_doubly_14a19470}, we show on the left the knot 14a19470 as the closure of $(\sigma_1 \sigma_2^{-1})^7$. After
flipping the braid inside the rectangle in the indicated direction, we get the diagram in the middle, and two symmetric Reidemeister 2 moves
yield the almost doubly symmetric diagram on the right (with the two exceptional crossings marked in gray).

\subsection{The number of distinct partial knots}
As mentioned before, Table \ref{tab:knotList} contains knot diagrams with different partial knots if these occurred as search results. 
For 59 knots we found diagrams with only one partial knot. For 33 knots there are two partial knots, for 11 knots there are three,
for 2 knots there are four, and there is one knot with five partial knots.
These results are experimental and we do not know if more partial knot variants exist in each of these cases.
The knot 12a427 has the highest number of partial knots and these are: $7_4$, $8_{21}$, $9_2$, $3_1 \sharp 4_1$, $3_1 \sharp 5_1$, all of them
having determinant 15. There are several interesting questions concerning partial knots and we note the following explicitly:

\begin{question}
Does a doubly symmetric strongly positive amphicheiral knot exist for which the partial knot in doubly symmetric diagrams is unique
but for which several partial knots occur in general symmetric union diagrams?
\end{question}

\subsection{Composite partial knots}
Doubly symmetric knots with composite partial knots occur frequently and the following templates generate them:
$t_2$, $t_3$, $t_{12}$, $t_{13}$, $t_{14}$, $t_{15}$, $t_{16}$. 
(Note, however, that in the templates $t_{13}$ and $t_{15}$ the knots 16n101996 and 16n102000 have the partial knot $3_1$ which is not composite. 
In this case the second factor is trivial and the partial knot should be read as $3_1 \sharp$ triv.) 

Two more remarks:
a) We did not find a doubly symmetric diagram representing a knot up to 16 crossings with a composite partial knot consisting of more than 2 prime factors.
b) Our notation for composite partial knots in the appendix is not completely precise: A partial knot written as $3_1 \sharp 3_1$ could either be $3_1 \sharp 3_1$ 
(with signature $\pm 4$) or $3_1 \sharp 3_1^*$, in which case the second factor is the mirror image of the first one and the total signature is 0.
The same ambiguity occurs for other partial knots, e.g. $3_1 \sharp 5_2$.

\subsection{Determinants and Alexander polynomials}
The determinants for the 131 prime strongly positive amphicheiral knots up to 16 crossings cover all squares of odd numbers between
1 and 53. The only knot with determinant $2809 = 53^2$ is 16a357272 and it has not yet been found to be doubly symmetric.
All other determinants between $1$ and $51^2$ occur for knots with doubly symmetric diagrams. The most frequent determinant is 441
with 9 knots (all of them doubly symmetric).

As the Alexander polynomial of strongly positive amphicheiral knots is always a square, we checked 
for doubly symmetric diagrams whether the square of the Alexander polynomial of one of the partial knots equals the 
Alexander polynomial of the knot. This is the case for roughly half of them (54 out of 106 knots). 

Examples: a) The Alexander polynomial of 14n22073 is $\Delta(t) = (1 - t^2 + t^4)^2$ and we found doubly symmetric diagrams for 14n22073 
with trivial partial knot only. The first knot in the knot table having Alexander polynomial $\Delta(t) = 1 - t^2 + t^4$ is 12n121 and 
it seems unlikely that 14n22073 can be a symmetric union with a partial knot with crossing number 12 or higher, as for instance 12n121.

b) In a similar example with a partial knot with crossing number 11 we find this relationship between the Alexander polynomials, though:
The doubly symmetric knot 16n869383 has Alexander polynomial $\Delta(t) = (2 - 4t + 3t^2 - 4t^3 + 2t^4)^2$ and we found 3 doubly symmetric 
diagrams, with partial knots $7_4$, $8_{21}$ and 11n79. The Alexander polynomial of 11n79 is $\Delta(t) = 2 - 4t + 3t^2 - 4t^3 + 2t^4$.

Remark: This relationship between the Alexander polynomials of a symmetric union and its partial knot always occurs for even twist numbers
(on the vertical axis). Doubly symmetric diagrams with even twist numbers are quite rare, though: Table \ref{tab:knotList} contains 
only the 8 cases 14n9732, 16n645918, 16n645926, 16n786382, 16n847920, 16n847983, 16n991381, 16n991505.

\begin{question}
Which topological properties distinguish doubly symmetric knots with even twist numbers from arbitrary doubly symmetric knots?
\end{question}

\begin{figure}[hbtp]
\centering
\includegraphics[scale=0.8]{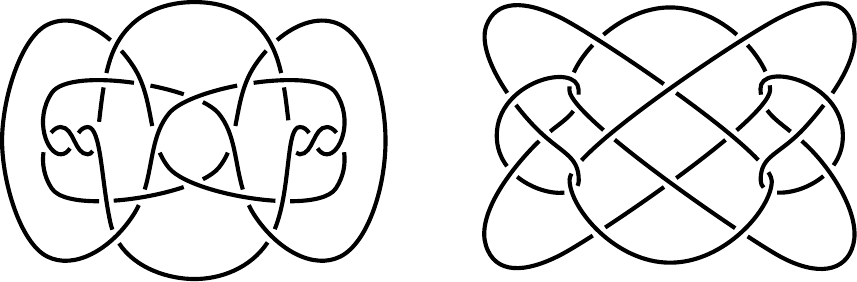}
\caption{Doubly symmetric knot diagrams for $16a_{107998}$, $t_{11}(0, 3, -1 \mid 1)$
and $16n_{845519}$, $t_{22}(0, 1, -1 \mid 1)$, both with 22 crossings}
\label{16a107998_16n845519}
\end{figure}

\subsection{Number of crossings in doubly symmetric diagrams}
We close the article with a comment on the crossing number of doubly symmetric diagrams. As mentioned before, alternating knots cannot 
realize their minimal crossing number in a doubly symmetric diagram. This is achieved, however, by the knot 14n9732 in $t_1(1, 0, \mid 2)$ and by 7
non-alternating knots with $c=16$: 16n428839, 16n451788, 16n645918, 16n645926, 16n847920, 16n991381 and 16n991505.
In the article \cite{Lamm2024} we will list exhaustively all doubly symmetric diagrams up to 18 crossings which are not reducible to less complicated ones
and we will formulate this as a theorem.

\begin{figure}[hbtp]
\centering
\includegraphics[scale=0.8]{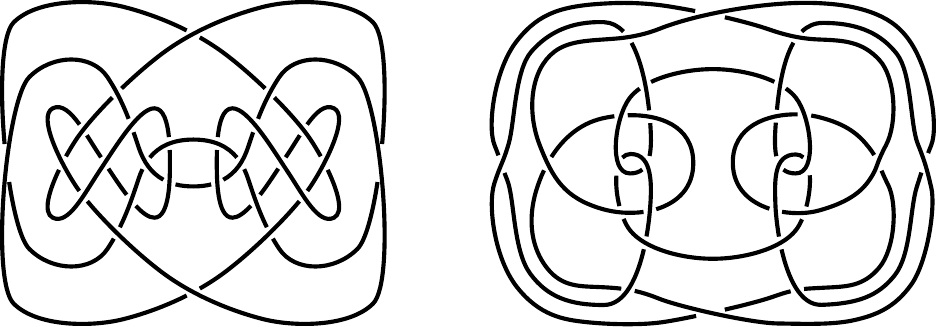}
\caption{Doubly symmetric knot diagrams for $16n_{225139}$, $t_{23}(-1, 1, -1, 1 \mid 1, 0)$ 
and $16n_{919894}$, $t_{21}(1, -1, -2, 0 \mid 1, 0)$, both with 26 crossings}
\label{16n225139_16n919894}
\end{figure}

The four doubly symmetric diagrams in Figures \ref{16a107998_16n845519} and \ref{16n225139_16n919894} additionally illustrate the opposite 
case that diagrams with quite a lot of crossings were found as minimal diagrams (22 in the first and 26 crossings in the second figure).
The maximal value of 28 crossings in a diagram in the appendix occurs only for one knot, and this is 16n925482 in the almost doubly
symmetric template $t_{28}$.

\section{Outlook}
The description of the systematic enumeration of doubly symmetric diagrams will be the topic of the article \cite{Lamm2024} 
(we decided not to include this in the present article because it would become too long). There, we will prove that no prime 
knots have doubly symmetric diagrams with $\le 12$ crossings in their symmetric diagram, and list the prime knots with doubly 
symmetric diagrams with 14 crossings (we find $10_{99}$, $10_{123}$, 12n706 and 14n9732) and with 16 and 18 crossings. 
We do that by enumerating all templates up to 4 crossings which are not reducible to less complicated ones; we use 
Knotscape \cite{Knotscape} to analyze the diagrams obtained by inserting crossings on the axes.

\section*{Acknowledgments}
I thank Marc Kegel and Steven Sivek for their help with my first steps with SnapPy, Nathan Dunfield for help concerning
the ribbon status of several knots and explanations concerning SnapPy, Brendan Owens for giving me a reference for the 
ribbon status of 14n19470, and Keegan Boyle for explaining the classification of symmetries for amphichei\-ral knots.

\small
\noindent
\begin{table}[p]
\begin{tabular}{|lrll|lrll|lrll|}
\hline 
$10_{99}$  &  81  & $6_1$             & $t_1$    & 16a20777  & 1521 & $3_1 \sharp 6_3$ & $t_{13}$ & 16n102009 &  441 & $3_1 \sharp 5_2$ & $t_{15}$ \\
           &      & $3_1 \sharp 3_1$  & $t_2$    & 16a51808  & 1369 & $9_{14}$         & $t_9$    & 16n102453 &  225 & $8_{21}$         & $t_7$    \\
           &      & $8_{20}$          & $t_5$    & 16a59437  & 2401 & $9_{27}$         & $t_9$    &           &      & $3_1 \sharp 5_1$ & $t_{14}$ \\
           &      & $9_1$             & $t_{17}$ & 16a66481  & 1225 & $9_{12}$         & $t_9$    & 16n106013 &  225 & $3_1 \sharp 5_1$ & $t_{16}$ \\
$10_{123}$ &  121 & $6_2$             & $t_1$    & 16a101002 & 1225 & $4_1 \sharp 5_2$ & $t_2$    & 16n204537 &   49 & $5_2$            & $t_7$    \\
           &      & $7_2$             & $t_9$    &           &      & $5_1 \sharp 5_2$ & $t_{12}$ &           &      & $7_1$            & $t_9$    \\ \cline{1-4}
12a427     &  225 & $7_4$             & $t_7$    & 16a107430 & 1089 & $8_{15}$         & $t_4$    & 16n225139 &  289 & $8_2$            & $t_{23}$ \\
           &      & $3_1 \sharp 4_1$  & $t_2$    &           &      & $9_{11}$         & $t_9$    & 16n225152 &  289 & $7_5$            & $t_9$    \\
           &      & $8_{21}$          & $t_4$    & 16a107998 & 1849 & $9_{22}$         & $t_{11}$ & 16n226358 &  625 & $8_8$            & $t_9$    \\
           &      & $3_1 \sharp 5_1$  & $t_{14}$ & 16a116332 & 1089 & $3_1 \sharp 6_2$ & $t_{12}$ & 16n226373 &   49 & $5_2$            & $t_9$    \\
           &      & $9_2$             & $t_{17}$ & 16a116809 & 1521 & $3_1 \sharp 6_3$ & $t_{12}$ &           &      & $7_1$            & $t_{10}$ \\ 
12a1019    &  361 & $7_6$             & $t_1$    &           &      & $9_{17}$         & $t_{21}$ & 16n268004 &  361 & $8_4$            & $t_9$    \\
12a1105    &  289 & $7_5$             & $t_1$    & 16a149574 &  729 & $3_1 \sharp 6_1$ & $t_2$    &           &      & $9_3$            & $t_9$    \\
           &      & $8_2$             & $t_9$    & 16a151960 & 1369 & $9_{36}$         & $t_{11}$ & 16n268599 &  121 & $7_2$            & $t_9$    \\
           &      & $10_1$            & $t_{17}$ & 16a152693 & 2025 & $9_{24}$         & $t_{11}$ & 16n272647 &  121 & $6_2$            & $t_{17}$ \\
12a1202    &  169 & $7_3$             & $t_1$    & 16a154813 & 1521 & $9_{17}$         & $t_{21}$ &           &      & $7_2$            & $t_9$    \\
           &      & $9_{43}$          & $t_5$    &           &      & $10_{144}$       & $t_4$    & 16n323632 &    9 & $3_1$            & $t_8$    \\ \cline{1-4}
12n706     &   49 & $5_2$             & $t_1$    & 16a168328 & 1681 & $10_{47}$        & $t_{19}$ &           &      & $8_{19}$         & $t_7$    \\
           &      & $7_1$             & $t_9$    & 16a193681 & 2209 & $9_{25}$         & $t_{11}$ & 16n428839 &  225 & $7_4$            & $t_7$    \\
           &      & $9_{42}$          & $t_4$    & 16a202161 &  625 & $4_1 \sharp 5_1$ & $t_{12}$ &           &      & $8_{21}$         & $t_4$    \\ \cline{1-4}
14a6002    &  625 & $4_1 \sharp 4_1$  & $t_2$    & 16a202163 & 1225 & $4_1 \sharp 5_2$ & $t_2$    &           &      & $9_2$            & $t_{17}$ \\
           &      & $4_1 \sharp 5_1$  & $t_{16}$ &           &      & $5_1 \sharp 5_2$ & $t_{14}$ &           &      & $10_{136}$       & $t_{22}$ \\
14a6398    &  529 & $8_7$             & $t_9$    & 16a259485 & 1089 & $8_{15}$         & $t_7$    & 16n451781 &  529 & $8_7$            & $t_{10}$ \\
14a8662    &  441 & $7_7$             & $t_4$    &           &      & $3_1 \sharp 6_2$ & $t_{14}$ & 16n451788 &   25 & $5_1$            & $t_1$    \\
           &      & $3_1 \sharp 5_2$  & $t_2$    & 16a259789 & 1089 & $3_1 \sharp 6_2$ & $t_{16}$ & 16n524438 &   49 & $9_{42}$         & $t_4$    \\
14a8872    &  729 & $8_{10}$          & $t_{11}$ & 16a293262 & 1089 & $3_1 \sharp 6_2$ & $t_{12}$ & 16n535891 &  225 & 11n13            & $t_{20}$ \\
           &      & $9_6$             & $t_{17}$ & 16a313024 & 2025 & $8_{18}$         & $t_4$    & 16n645918 &   81 & $3_1 \sharp 3_1$ & $t_2$    \\
           &      & $9_{48}$          & $t_4$    &           &      & $10_{164}$       & $t_4$    & 16n645926 &   81 & $3_1 \sharp 3_1$ & $t_3$    \\
14a16309   &  625 & $9_{49}$          & $t_4$    & 16a313458 &  729 & $3_1 \sharp 6_1$ & $t_2$    & 16n780729 &  529 & $8_6$            & $t_{10}$ \\
14a16311   &  289 & $8_3$             & $t_1$    & 16a314171 & 1225 & $8_{16}$         & $t_4$    &           &      & $10_2$           & $t_{17}$ \\
           &      & $10_{130}$        & $t_5$    &           &      & $9_{12}$         & $t_1$    & 16n780731 &    1 & triv             & $t_{10}$ \\
14a16437   &  729 & $9_6$             & $t_9$    & 16a326803 & 1521 & $10_{165}$       & $t_4$    & 16n786382 &    1 & triv             & $t_8$    \\
14a17173   &  841 & $8_{12}$          & $t_1$    & 16a329540 & 2025 & $10_{158}$       & $t_4$    & 16n797553 &  361 & $8_4$            & $t_9$    \\
14a18187   &  729 & $8_{11}$          & $t_1$    & 16a330218 & 1089 & $8_{15}$         & $t_4$    & 16n847920 &   49 & $5_2$            & $t_1$    \\
           &      & $9_{47}$          & $t_4$    & 16a332442 & 1369 & $9_{36}$         & $t_{11}$ & 16n847983 &   25 & $4_1$            & $t_9$    \\
14a18362   &  361 & $8_4$             & $t_9$    & 16a332444 & 1521 & $9_{15}$         & $t_1$    &           &      & $5_1$            & $t_8$    \\
14a18676   &  225 & $7_4$             & $t_4$    & 16a340727 &  961 & $10_{46}$        & $t_{19}$ & 16n847988 &  289 & $7_5$            & $t_9$    \\
           &      & $3_1 \sharp 5_1$  & $t_{12}$ & 16a340770 & 1369 & $10_6$           & $t_9$    & 16n849519 &  441 & $7_7$            & $t_{22}$ \\
14a18680   &  441 & $7_7$             & $t_7$    &           &      & 11a234           & $t_{17}$ & 16n855304 &  289 & $9_{44}$         & $t_5$    \\
           &      & $3_1 \sharp 5_2$  & $t_2$    & 16a354511 & 1225 & $8_{16}$         & $t_4$    & 16n858257 &  441 & $7_7$            & $t_1$    \\
14a18723   &  529 & $8_6$             & $t_1$    &        &         & $10_{162}$       & $t_4$    &           &      & $8_5$            & $t_7$    \\
           &      & $9_5$             & $t_{10}$ & 16a354563 & 1681 & $9_{18}$         & $t_1$    & 16n868471 &  169 & $7_3$            & $t_{10}$ \\
           &      & $10_2$            & $t_{17}$ & 16a356133 & 2601 & $10_{163}$       & $t_4$    &           &      & $8_1$            & $t_9$    \\
14a19472   &  441 & $7_7$             & $t_4$    & 16a356811 &  961 & $9_9$            & $t_9$    & 16n869383 &  225 & $7_4$            & $t_1$    \\
           &      & $8_5$             & $t_{11}$ & 16a356843 &  441 & $8_5$            & $t_4$    &           &      & $8_{21}$         & $t_4$    \\
           &      & $9_4$             & $t_{18}$ &           &      & $9_4$            & $t_1$	  &           &      & 11n79            & $t_4$    \\ \cline{1-4}
14n8213    &  169 & $6_3$             & $t_1$    &           &      & 11n64            & $t_5$    & 16n872167 &   81 & $8_{20}$         & $t_6$    \\ 
           &      & $7_3$             & $t_{10}$ & 16a358693 & 1681 & $9_{20}$         & $t_9$    &           &      & $9_1$            & $t_{17}$ \\
           &      & $8_1$             & $t_9$    & 16a368293 & 1521 & $9_{16}$         & $t_{11}$ &           &      & $9_{46}$         & $t_4$    \\ 
14n9732    &   25 & $4_1$             & $t_7$    &           &      & $10_9$           & $t_{18}$ & 16n872172 &  361 & $7_6$            & $t_1$    \\ 
           &      & $5_1$             & $t_1$    & 16a371892 &  841 & $9_7$            & $t_1$    &           &      & $9_3$            & $t_9$    \\ 
           &      & $10_{132}$        & $t_{21}$ &           &      & $10_8$           & $t_9$    & 16n919894 &  441 & 11n118           & $t_{21}$ \\ \cline{5-8}
14n22073   &    1 & triv              & $t_8$    & 16n101996 &    9 & $3_1$            & $t_{13}$ & 16n925976 &    9 & $8_{19}$         & $t_7$    \\ 
14n25903   &  121 & $6_2$             & $t_1$    &           &      & $8_{19}$         & $t_7$    & 16n958050 &  625 & $8_9$            & $t_{10}$ \\
           &      & $10_{125}$        & $t_{20}$ & 16n102000 &    9 & $3_1$            & $t_{15}$ & 16n988939 &    1 & triv             & $t_{18}$ \\
	       &      & $10_{128}$        & $t_4$    & 16n102007 &  441 & $8_5$            & $t_7$    & 16n991381 &   81 & $6_1$            & $t_1$    \\ \cline{1-4}
16a19960   & 1521 & $3_1 \sharp 6_3$  & $t_{15}$ &           &      & $3_1 \sharp 5_2$ & $t_{13}$ & 16n991505 &  121 & $6_2$            & $t_1$    \\
\hline 
\end{tabular}
\caption{A list of 171 diagrams for the 106 prime strongly positive amphicheiral knots for which doubly symmetric diagrams have been found. 
The second and third columns contain determinant and partial knot and the fourth column indicates the template with the diagram information.}
\label{tab:knotList}
\end{table}

\normalsize

\small
\noindent
\begin{table}[htbp]
\begin{tabular}{|lrl|lrl|lrl|}
\hline 
14a19470  &  841 & $t_{24}$ & 16a288139 & 1681 & $t_{24}$ & 16a339566 & 2209 & $t_{26}$ \\
14a19517  &  529 & $t_{28}$ & 16a289378 & 2209 & $t_{24}$ & 16a354193 & 2025 & $t_{28}$ \\ \cline{1-3}
16a66232  & 1849 & $t_{25}$ & 16a309401 & 2401 & $t_{28}$ & 16a357272 & 2809 & $t_{26}$ \\
16a106900 & 1369 & $t_{27}$ & 16a312423 & 1521 & $t_{25}$ & 16a359219 & 1225 & $t_{27}$ \\ \cline{7-9}
16a150121 & 1681 & $t_{25}$ & 16a324978 & 1089 & $t_{25}$ & 16n404143 &  225 & $t_{25}$ \\
16a155718 &  841 & $t_{27}$ & 16a326826 & 1849 & $t_{26}$ & 16n526315 &  169 & $t_{25}$ \\
16a156092 &  961 & $t_{25}$ & 16a329518 & 1849 & $t_{25}$ & 16n918713 &  361 & $t_{25}$ \\
16a275557 & 1369 & $t_{25}$ & 16a330763 & 2601 & $t_{26}$ & 16n925482 &   25 & $t_{28}$ \\
16a275642 & 2025 & $t_{25}$ &           &      &          &           &      & \\
\hline 
\end{tabular}
\caption{A list of the 25 prime strongly positive amphicheiral knots for which doubly symmetric diagrams have not been found. 
The second column contains the determinant and the third column indicates the (almost doubly symmetric) template with the diagram information.}
\label{tab:knotList2}
\end{table}

\normalsize

\clearpage
\section{Appendix}
The Appendix contains the templates generating the 196 (171 + 25) diagrams of strongly positive amphicheiral prime knots
with crossing number $\le 16$. The twist numbers correspond to the dotted lines in the templates read from left to right (first part) 
and from top to bottom (second part). In the second column the partial knot is indicated (notation $\partial$). 
The first four pages contain the doubly symmetric diagrams and the fifth page shows the almost doubly symmetric diagrams.

\small

\vspace{1.0cm}
\hspace{-1cm}
\noindent
\parbox[b]{4.2cm}{
\centering
\mbox{} \\
\includegraphics[scale=\scaling]{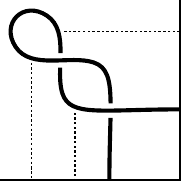} \\
template $t_1$ \\
}
\quad
\parbox[b]{4.2cm}{
\centering
\mbox{} \\
\includegraphics[scale=\scaling]{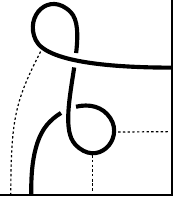} \\
template $t_2$ \\
}
\quad
\parbox[b]{4.2cm}{
\centering
\mbox{} \\
\includegraphics[scale=\scaling]{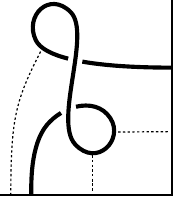} \\
template $t_3$ \\
}

\noindent
\parbox[t]{4.2cm}{
\centering
\mbox{} \\
\begin{tabular}{l@{, $\partial$}l@{\;=\;}r@{, }r@{ $\mid$ }r}
$10_{99}$      & $6_1$    &( 2&  0& 1) \\
$10_{123}$     & $6_2$    &( 1&  1& 1) \\
$12a_{1019}$   & $7_6$    &( 2&  1& 1) \\
$12a_{1105}$   & $7_5$    &( 1&  2& 1) \\
$12a_{1202}$   & $7_3$    &( 3&  0& 1) \\
$12n_{706}$    & $5_2$    &(-2&  0& 1) \\
$14a_{16311}$  & $8_3$    &( 4&  0& 1) \\
$14a_{17173}$  & $8_{12}$ &( 2&  2& 1) \\
$14a_{18187}$  & $8_{11}$ &( 3&  1& 1) \\
$14a_{18723}$  & $8_6$    &( 1&  3& 1) \\
$14n_{8213}$   & $6_3$    &(-2&  1& 1) \\
$14n_{9732}$   & $5_1$    &( 1&  0& 2) \\
$14n_{25903}$  & $6_2$    &(-3&  0& 1) \\
$16a_{314171}$ & $9_{12}$ &( 4&  1& 1) \\
$16a_{332444}$ & $9_{15}$ &( 2&  3& 1) \\
$16a_{354563}$ & $9_{18}$ &( 3&  2& 1) \\
$16a_{356843}$ & $9_4$    &( 5&  0& 1) \\
$16a_{371892}$ & $9_7$    &( 1&  4& 1) \\
$16n_{451788}$ & $5_1$    &( 1&  0& 3) \\
$16n_{847920}$ & $5_2$    &(-2&  0& 2) \\
$16n_{858257}$ & $7_7$    &(-3&  1& 1) \\
$16n_{869383}$ & $7_4$    &(-4&  0& 1) \\
$16n_{872172}$ & $7_6$    &(-2&  2& 1) \\
$16n_{991381}$ & $6_1$    &( 2&  0& 2) \\
$16n_{991505}$ & $6_2$    &( 1&  1& 2) \\
\end{tabular}
}
\quad
\noindent
\parbox[t]{4.2cm}{
\centering
\mbox{} \\
\begin{tabular}{l@{, $\partial$}l@{\;=\;}r@{, }r@{ $\mid$ }r}
$10_{99}$       & $3_1 \sharp 3_1$  &( 1& -1& 1) \\
$12a_{427}$     & $3_1 \sharp 4_1$  &( 1& -2& 1) \\
$14a_{6002}$    & $4_1 \sharp 4_1$  &( 2& -2& 1) \\
$14a_{8662}$    & $3_1 \sharp 5_2$  &( 1& -3& 1) \\
$14a_{18680}$   & $3_1 \sharp 5_2$  &( 3& -1& 1) \\
$16a_{101002}$  & $4_1 \sharp 5_2$  &( 2& -3& 1) \\
$16a_{149574}$  & $3_1 \sharp 6_1$  &( 1& -4& 1) \\
$16a_{202163}$  & $4_1 \sharp 5_2$  &( 3& -2& 1) \\
$16a_{313458}$  & $3_1 \sharp 6_1$  &( 4& -1& 1) \\
$16n_{645918}$  & $3_1 \sharp 3_1$  &( 1& -1& 2) \\
\end{tabular}
}
\quad
\parbox[t]{4.2cm}{
\centering
\mbox{} \\
\vspace{-0.1cm}
\begin{tabular}{l@{, $\partial$}l@{\;=\;}r@{, }r@{ $\mid$ }r}
$16n_{645926}$  & $3_1 \sharp 3_1$  &(-1& -1& 2) \\
\end{tabular}
}

\vspace{-2.3cm}
\hspace{4cm}
\noindent
\parbox[t]{4.2cm}{
\centering
\mbox{} \\
\includegraphics[scale=\scaling]{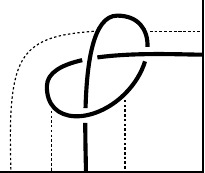} \\
template $t_4$ \\
\centering
\mbox{} \\
\begin{tabular}{l@{, $\partial$}l@{\;=\;}r@{, }r@{, }r@{ $\mid$ }r}
$12a_{427}$      & $8_{21}$    &(-2&  0&  0& 1) \\
$12n_{706}$      & $9_{42}$    &( 0&  1&  2& 1) \\
$14a_{8662}$     & $7_7$       &( 1&  0&  0& 1) \\
$14a_{8872}$     & $9_{48}$    &(-3&  0&  0& 1) \\
$14a_{16309}$    & $9_{49}$    &(-2&  0& -1& 1) \\
$14a_{18187}$    & $9_{47}$    &(-2&  1&  0& 1) \\
$14a_{18676}$    & $7_4$       &( 0&  0& -1& 1) \\
$14a_{19472}$    & $7_7$       &( 0&  1&  0& 1) \\
$14n_{25903}$    & $10_{128}$  &( 0&  2&  2& 1) \\
\end{tabular}
}
\quad
\parbox[t]{4.2cm}{
\centering
\mbox{} \\
\vspace{0.07cm}
\begin{tabular}{l@{, $\partial$}l@{\;=\;}r@{, }r@{, }r@{ $\mid$ }r}
$16a_{107430}$   & $8_{15}$    &( 2&  0&  0& 1) \\
$16a_{154813}$   & $10_{144}$  &(-4&  0&  0& 1) \\
$16a_{313024}$   & $8_{18}$    &( 1&  1&  0& 1) \\
$16a_{313024}$   & $10_{164}$  &(-2&  1& -1& 1) \\
$16a_{314171}$   & $8_{16}$    &( 0&  1& -1& 1) \\
$16a_{326803}$   & $10_{165}$  &(-2&  2&  0& 1) \\
$16a_{329540}$   & $10_{158}$  &(-3&  0& -1& 1) \\
$16a_{330218}$   & $8_{15}$    &( 0&  2&  0& 1) \\
$16a_{354511}$   & $8_{16}$    &( 1&  0& -1& 1) \\
$16a_{354511}$   & $10_{162}$  &(-2&  0& -2& 1) \\
$16a_{356133}$   & $10_{163}$  &(-3&  1&  0& 1) \\
$16a_{356843}$   & $8_5$       &( 0&  0& -2& 1) \\
$16n_{428839}$   & $8_{21}$    &( 1&  1&  1& 1) \\
$16n_{524438}$   & $9_{42}$    &( 1&  0&  2& 1) \\
$16n_{869383}$   & $8_{21}$    &( 0& -2&  0& 1) \\
$16n_{869383}$   & $11n_{79}$  &( 0&  3&  2& 1) \\
$16n_{872167}$   & $9_{46}$    &( 0&  0&  3& 1) \\
\end{tabular}
}

\hspace{-1cm}
\noindent
\parbox[b]{4.5cm}{
\centering
\mbox{} \\
\includegraphics[scale=\scaling]{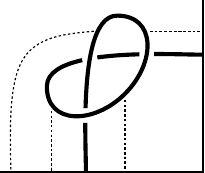} \\
template $t_5$ \\
}
\hspace{1.2cm}
\parbox[b]{4.5cm}{
\centering
\mbox{} \\
\includegraphics[scale=\scaling]{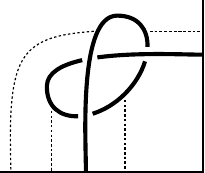} \\
template $t_6$ \\
}

\noindent
\parbox[t]{4.5cm}{
\centering
\mbox{} \\
\begin{tabular}{l@{, $\partial$}l@{\;=\;}r@{, }r@{, }r@{ $\mid$ }r}
$10_{99}$      & $8_{20}$    &(-1&  1&  0& 1) \\
$12a_{1202}$   & $9_{43}$    &(-1&  2&  0& 1) \\
$14a_{16311}$  & $10_{130}$  &(-1&  3&  0& 1) \\
$16a_{356843}$ & $11n_{64}$  &(-1&  4&  0& 1) \\
$16n_{855304}$ & $9_{44}$    &(-2&  1&  0& 1) \\
\end{tabular}
}
\hspace{1.2cm}
\noindent
\parbox[t]{4.5cm}{
\centering
\mbox{} \\
\begin{tabular}{l@{, $\partial$}l@{\;=\;}r@{, }r@{, }r@{ $\mid$ }r}
$16n_{872167}$      & $8_{20}$ &( 1& -1&  0& 1) \\
\end{tabular}
}

\vspace{1.3cm}
\hspace{-1cm}
\parbox[t]{3.5cm}{
\centering
\mbox{} \\
\includegraphics[scale=\scaling]{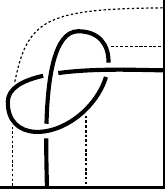} \\
template $t_7$ \\
}
\parbox[t]{4.5cm}{
\centering
\mbox{} \\
\begin{tabular}{l@{, $\partial$}l@{\;=\;}r@{, }r@{ $\mid$ }r@{, }r}
$12a_{427}$    & $7_4$    &( 0& -1&-1& 1) \\
$14a_{18680}$  & $7_7$    &( 1&  0&-1& 1) \\
$14n_{9732}$   & $4_1$    &(-1& -1&-1& 0) \\
$16a_{259485}$ & $8_{15}$ &( 2&  0&-1& 1) \\
$16n_{101996}$ & $8_{19}$ &( 0&  2&-1& 1) \\
$16n_{102007}$ & $8_5$    &( 0& -2&-1& 1) \\
$16n_{102453}$ & $8_{21}$ &(-2&  0&-1& 1) \\
$16n_{204537}$ & $5_2$    &(-1& -2&-1& 0) \\
$16n_{323632}$ & $8_{19}$ &( 0&  2&-1& 0) \\
$16n_{428839}$ & $7_4$    &( 0& -1&-1& 0) \\
$16n_{858257}$ & $8_5$    &( 0& -2&-1& 2) \\
$16n_{925976}$ & $8_{19}$ &( 0&  2&-1& 2) \\
\end{tabular}
}
\hspace{1.2cm}
\parbox[t]{4.5cm}{
\centering
\mbox{} \\
\includegraphics[scale=\scaling]{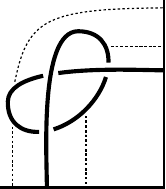} \\
template $t_8$ \\
\mbox{} \\
\vspace{0.1cm}
\begin{tabular}{l@{, $\partial$}l@{\;=\;}r@{, }r@{ $\mid$ }r@{, }r}
$14n_{22073}$  & triv  &( 0&  1& 1& 0) \\
$16n_{323632}$ & $3_1$ &( 0&  2& 1& 0) \\
$16n_{786382}$ & triv  &( 0&  1& 2& 0) \\
$16n_{847983}$ & $5_1$ &(-1&  1& 1& 0) \\
\end{tabular}
}

\vspace{1.3cm}
\hspace{-1cm}
\parbox[t]{4.5cm}{
\centering
\mbox{} \\
\includegraphics[scale=\scaling]{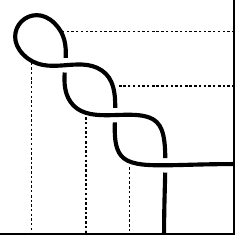} \\
template $t_9$ \\
\mbox{} \\
\begin{tabular}{l@{, $\partial$}l@{\;=\;}r@{, }r@{, }r@{ $\mid$ }r@{, }r}
$10_{123}$     & $7_2$    &(-2&  0&  0& 0&  1) \\
$12a_{1105}$   & $8_2$    &( 1&  0&  1& 1&  0) \\
$12n_{706}$    & $7_1$    &( 1&  0&  0& 0&  1) \\
$14a_{6398}$   & $8_7$    &(-2&  1&  0& 0&  1) \\
$14a_{16437}$  & $9_6$    &( 1&  0&  2& 1&  0) \\
$14a_{18362}$  & $8_4$    &( 1&  1&  0& 0&  1) \\
$14n_{8213}$   & $8_1$    &( 2&  0&  0& 0&  1) \\
$16a_{51808}$  & $9_{14}$ &(-3&  1&  0& 0&  1) \\
$16a_{59437}$  & $9_{27}$ &(-2&  1&  1& 0&  1) \\
$16a_{66481}$  & $9_{12}$ &(-2&  2&  0& 0&  1) \\
\end{tabular}
}
\quad
\parbox[t]{4.5cm}{
\centering
\mbox{} \\
\vspace{0.9cm}
\begin{tabular}{l@{, $\partial$}l@{\;=\;}r@{, }r@{, }r@{ $\mid$ }r@{, }r}
$16a_{107430}$ & $9_{11}$ &( 2&  1&  0& 0&  1) \\
$16a_{340770}$ & $10_6$   &( 1&  0&  3& 1&  0) \\
$16a_{356811}$ & $9_9$    &( 1&  2&  0& 0&  1) \\
$16a_{358693}$ & $9_{20}$ &( 1&  1&  1& 0&  1) \\
$16a_{371892}$ & $10_8$   &(-5&  0&  0& 0&  1) \\
$16n_{204537}$ & $7_1$    &( 1&  0&  0&-2&  1) \\
$16n_{225152}$ & $7_5$    &( 1& -2&  0& 0&  1) \\
$16n_{226358}$ & $8_8$    &(-2&  0&  1& 0&  1) \\
$16n_{226373}$ & $5_2$    &( 1& -1&  1& 0&  1) \\
$16n_{268004}$ & $8_4$    &( 1&  1&  0&-1&  1) \\
$16n_{268004}$ & $9_3$    &( 3&  0&  0&-1&  1) \\
$16n_{268599}$ & $7_2$    &(-2&  0&  0& 1&  0) \\
$16n_{272647}$ & $7_2$    &(-2&  0&  0&-1&  1) \\
$16n_{797553}$ & $8_4$    &( 1&  1&  0& 1&  0) \\
$16n_{847983}$ & $4_1$    &( 1& -1&  0& 0&  2) \\
$16n_{847988}$ & $7_5$    &( 1& -2&  0&-1&  1) \\
$16n_{868471}$ & $8_1$    &( 2&  0&  0& 1&  0) \\
$16n_{872172}$ & $9_3$    &( 3&  0&  0& 0&  1) \\
\end{tabular}
}
\hspace{1.0cm}
\parbox[t]{5cm}{
\centering
\mbox{} \\
\includegraphics[scale=\scaling]{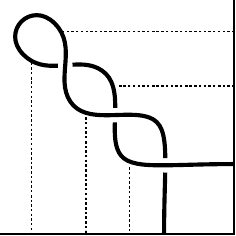} \\
template $t_{10}$ \\
\mbox{} \\
\begin{tabular}{l@{, $\partial$}l@{\;=\;}r@{, }r@{, }r@{ $\mid$ }r@{, }r}
$14a_{18723}$   & $9_5$    &(-3&  1&  0& 0&  1) \\
$14n_{8213}$    & $7_3$    &(-1& -1&  0& 1&  0) \\
$16n_{226373}$  & $7_1$    &(-1& -1& -2& 1&  1) \\
$16n_{451781}$  & $8_7$    &(-1& -1&  1& 1&  0) \\
$16n_{780729}$  & $8_6$    &(-2& -1&  0& 1&  0) \\
$16n_{780731}$  & triv     &(-1&  0&  1& 1&  1) \\
$16n_{868471}$  & $7_3$    &(-1& -1&  0& 0&  1) \\
$16n_{958050}$  & $8_9$    &(-1& -2&  0& 1&  0) \\
\end{tabular}
}

\noindent
\parbox[t]{6.5cm}{
\centering
\mbox{} \\
\includegraphics[scale=\scaling]{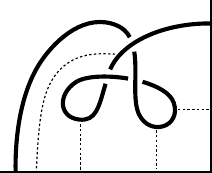} \\
template $t_{11}$ \\
}
\hspace{-0.9cm}
\parbox[t]{4.2cm}{
\centering
\mbox{} \\
\begin{tabular}{l@{, $\partial$}l@{\;=\;}r@{, }r@{, }r@{ $\mid$ }r}
$14a_{8872}$      & $8_{10}$    &( 0&  2& -1& 1) \\
$14a_{19472}$     & $8_5$       &( 0& -1& -1& 1) \\
$16a_{107998}$    & $9_{22}$    &( 0&  3& -1& 1) \\
$16a_{151960}$    & $9_{36}$    &( 0& -2& -1& 1) \\
$16a_{152693}$    & $9_{24}$    &(-1&  2& -1& 1) \\
$16a_{193681}$    & $9_{25}$    &( 0&  2& -2& 1) \\
$16a_{332442}$    & $9_{36}$    &( 0& -1& -2& 1) \\
$16a_{368293}$    & $9_{16}$    &(-1& -1& -1& 1) \\
\end{tabular}
}

\vspace{1.3cm}
\parbox[t]{5.8cm}{
\centering
\mbox{} \\
\includegraphics[scale=\scaling]{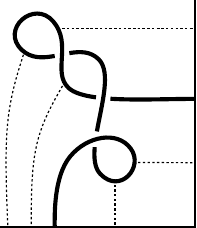} \\
template $t_{12}$ \\
}
\hspace{-0.7cm}
\parbox[t]{4.5cm}{
\centering
\mbox{} \\
\begin{tabular}{l@{, $\partial$}l@{\;=\;}r@{, }r@{, }r@{ $\mid$ }r@{, }r}
$14a_{18676}$    & $3_1 \sharp 5_1$ &(-1&  0&  1& 0&  1) \\
$16a_{101002}$   & $5_1 \sharp 5_2$ &(-1&  0&  3&-1&  1) \\
$16a_{116332}$   & $3_1 \sharp 6_2$ &( 3&  0&  1& 0&  1) \\
$16a_{116809}$   & $3_1 \sharp 6_3$ &( 2& -1&  1& 0&  1) \\
$16a_{202161}$   & $4_1 \sharp 5_1$ &(-1&  0&  2& 0&  1) \\
$16a_{293262}$   & $3_1 \sharp 6_2$ &(-1& -1&  1& 0&  1) \\
\end{tabular}
}

\vspace{1.3cm}
\parbox[t]{5.5cm}{
\centering
\mbox{} \\
\includegraphics[scale=\scaling]{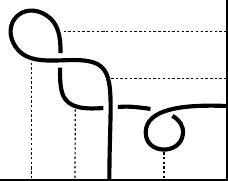} \\
template $t_{13}$ \\
\mbox{} \\
\vspace{0.1cm}
\begin{tabular}{l@{, $\partial$}l@{\;=\;}r@{, }r@{, }r@{ $\mid$ }r@{, }r}
$16a_{20777}$  & $3_1 \sharp 6_3$ &( 1& -2& -1& 0& 1) \\
$16n_{101996}$ & $3_1$            &( 1&  0& -1& 0& 1) \\
$16n_{102007}$ & $3_1 \sharp 5_2$ &( 1& -1& -1& 0& 1) \\
\end{tabular}
}
\quad
\parbox[t]{5.5cm}{
\centering
\mbox{} \\
\includegraphics[scale=\scaling]{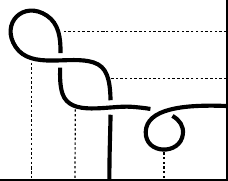} \\
template $t_{14}$ \\
\mbox{} \\
\vspace{0.1cm}
\begin{tabular}{l@{, $\partial$}l@{\;=\;}r@{, }r@{, }r@{ $\mid$ }r@{, }r}
$12a_{427}$    & $3_1 \sharp 5_1$ &( 1&  0& -1& 1& 0) \\
$16a_{202163}$ & $5_1 \sharp 5_2$ &( 1&  0& -3& 1& 0) \\
$16a_{259485}$ & $3_1 \sharp 6_2$ &( 1&  1& -1& 0& 1) \\
$16n_{102453}$ & $3_1 \sharp 5_1$ &( 1&  0& -1& 0& 1) \\
\end{tabular}
}

\vspace{1.3cm}
\parbox[t]{5.5cm}{
\centering
\mbox{} \\
\includegraphics[scale=\scaling]{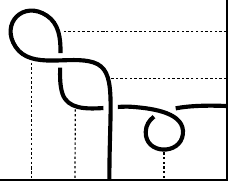} \\
template $t_{15}$ \\
\mbox{} \\
\vspace{0.1cm}
\begin{tabular}{l@{, $\partial$}l@{\;=\;}r@{, }r@{, }r@{ $\mid$ }r@{, }r}
$16a_{19960}$    & $3_1 \sharp 6_3$ &( 1& -2&  1& 0& 1) \\
$16n_{102000}$   & $3_1$            &( 1&  0&  1& 0& 1) \\
$16n_{102009}$   & $3_1 \sharp 5_2$ &( 1& -1&  1& 0& 1) \\
\end{tabular}
}
\quad
\parbox[t]{5.5cm}{
\centering
\mbox{} \\
\includegraphics[scale=\scaling]{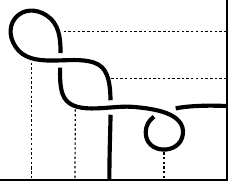} \\
template $t_{16}$ \\
\mbox{} \\
\vspace{0.1cm}
\begin{tabular}{l@{, $\partial$}l@{\;=\;}r@{, }r@{, }r@{ $\mid$ }r@{, }r}
$14a_{6002}$   & $4_1 \sharp 5_1$ &( 1&  0&  2& 1& 0) \\
$16a_{259789}$ & $3_1 \sharp 6_2$ &( 1&  1&  1& 0& 1) \\
$16n_{106013}$ & $3_1 \sharp 5_1$ &( 1&  0&  1& 0& 1) \\
\end{tabular}
}

\hspace{-1.2cm}
\parbox[t]{4.2cm}{
\centering
\mbox{} \\
\includegraphics[scale=\scaling]{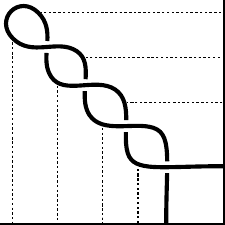} \\
template $t_{17}$ \\
}
\parbox[t]{4.5cm}{
\mbox{} \\
\begin{tabular}{l@{, $\partial$}l@{\;=\;}r@{, }r@{, }r@{, }r@{ $\mid$ }r@{, }r@{, }r}
$10_{99}$      & $9_1$       &( 1& 0& 0& 0& 0& -1& 1) \\
$12a_{1105}$   & $10_1$      &( 2& 0& 0& 0& 0& -1& 1) \\
$12a_{427}$    & $9_2$       &(-2& 0& 0& 0& 0& -1& 1) \\
$14a_{8872}$   & $9_6$       &( 1&-2& 0& 0& 0&  0& 1) \\
$14a_{18723}$  & $10_2$      &( 1& 0& 0& 1& 1&  0& 0) \\
$16a_{340770}$ & $11a_{234}$ &( 1& 0& 0& 2& 1&  0& 0) \\
$16n_{272647}$ & $6_2$       &( 1& 0&-1& 0& 1&  0& 1) \\
$16n_{428839}$ & $9_2$       &(-2& 0& 0& 0& 0&  0& 1) \\
$16n_{780729}$ & $10_2$      &( 1& 0& 0& 1& 1& -1& 0) \\
$16n_{872167}$ & $9_1$       &( 1& 0& 0& 0& 0&  0& 1) \\
\end{tabular}
}
\hspace{1.4cm}
\parbox[t]{4.5cm}{
\centering
\mbox{} \\
\includegraphics[scale=\scaling]{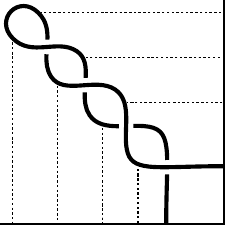} \\
template $t_{18}$ \\
\mbox{} \\
\vspace{0.3cm}
\begin{tabular}{l@{, $\partial$}l@{\;=\;}r@{, }r@{, }r@{, }r@{ $\mid$ }r@{, }r@{, }r}
$14a_{19472}$    & $9_4$  &( 1& 0& -1& -1& 0& 0& 1) \\
$16a_{368293}$   & $10_9$ &( 1& 1& -1& -1& 0& 0& 1) \\
$16n_{988939}$   & triv   &( 1& 0&  0& -1& 1& 0& 1) \\
\end{tabular}
}

\vspace{1cm}
\parbox[t]{5.5cm}{
\centering
\mbox{} \\
\includegraphics[scale=\scaling]{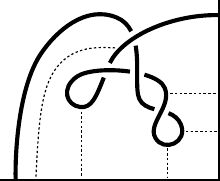} \\
template $t_{19}$ \\
\mbox{} \\
\vspace{0.1cm}
\begin{tabular}{l@{, $\partial$}l@{\;=\;}r@{, }r@{, }r@{ $\mid$ }r@{, }r}
$16a_{168328}$  & $10_{47}$ &( 0&  2& -1& 0& 1) \\
$16a_{340727}$  & $10_{46}$ &( 0& -1& -1& 0& 1) \\
\end{tabular}
}
\quad
\parbox[t]{5.5cm}{
\centering
\mbox{} \\
\includegraphics[scale=\scaling]{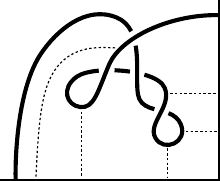} \\
template $t_{20}$ \\
\mbox{} \\
\vspace{0.1cm}
\begin{tabular}{l@{, $\partial$}l@{\;=\;}r@{, }r@{, }r@{ $\mid$ }r@{, }r}
$14n_{25903}$   & $10_{125}$ &( 0& 1& -1& 0& 1) \\
$16n_{535891}$  & $11n_{13}$ &( 0& 2& -1& 0& 1) \\
\end{tabular}
}

\vspace{1cm}
\parbox[t]{5.5cm}{
\centering
\mbox{} \\
\includegraphics[scale=\scaling]{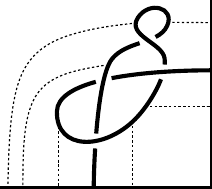} \\
template $t_{21}$ \\
\mbox{} \\
\vspace{0.1cm}
\begin{tabular}{l@{, $\partial$}l@{\;=\;}r@{, }r@{, }r@{, }r@{ $\mid$ }r@{, }r}
$14n_{9732}$    & $10_{132}$  &( 1& -1&  0&-1& 1& 0) \\
$16a_{116809}$  & $9_{17}$    &( 1&  0&  0& 0& 1& 0) \\
$16a_{154813}$  & $9_{17}$    &( 1&  0&  0& 0& 1&-1) \\
$16n_{919894}$  & $11n_{118}$ &( 1& -1& -2& 0& 1& 0) \\
\end{tabular}
}
\quad
\parbox[t]{5.5cm}{
\centering
\mbox{} \\
\includegraphics[scale=\scaling]{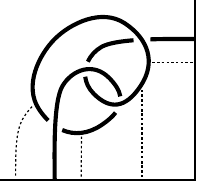} \\
template $t_{22}$ \\
\mbox{} \\
\vspace{0.1cm}
\begin{tabular}{l@{, $\partial$}l@{\;=\;}r@{, }r@{, }r@{ $\mid$ }r}
$16n_{428839}$  & $10_{136}$ &(0 &-1& -1& 1) \\
$16n_{849519}$  & $7_7$      &(0 & 1& -1& 1) \\
\end{tabular}
}

\enlargethispage{0.5cm}

\vspace{1cm}
\parbox[t]{5.5cm}{
\centering
\mbox{} \\
\includegraphics[scale=\scaling]{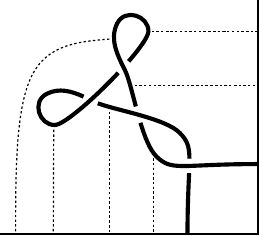} \\
template $t_{23}$ \\
\mbox{} \\
\vspace{0.1cm}
\begin{tabular}{l@{, $\partial$}l@{\;=\;}r@{, }r@{, }r@{, }r@{ $\mid$ }r@{, }r}
$16n_{225139}$  & $8_2$ &(-1& 1& -1& 1& 1& 0) \\
\end{tabular}
}

\parbox[t]{5cm}{
\centering
\mbox{} \\
\includegraphics[scale=\scaling]{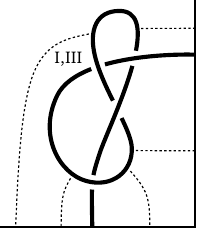} \\
template $t_{24}$ \\
\mbox{} \\
\begin{tabular}{l@{\;=\;}r@{, }r@{, }r@{ $\mid$ }r@{, }r}
$14a_{19470}$  &( 0&  0& -1&  1&  0) \\
$16a_{288139}$ &(-2&  1& -1&  0& -1) \\
$16a_{289378}$ &(-1&  1& -2&  0& -1) \\
\end{tabular}
}
\quad
\parbox[t]{5cm}{
\centering
\mbox{} \\
\includegraphics[scale=\scaling]{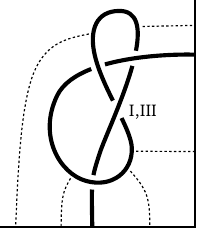} \\
template $t_{25}$ \\
\mbox{} \\
\begin{tabular}{l@{\;=\;}r@{, }r@{, }r@{ $\mid$ }r@{, }r}
$16a_{66232}$  &(  0& -2&  0&  0& -1) \\
$16a_{150121}$ &(  0&  1&  0&  0& -1) \\
$16a_{156092}$ &(  0&  2&  0& -1&  0) \\
$16a_{275557}$ &(  1& -2&  0&  1&  0) \\
$16a_{275642}$ &(  0& -2&  0&  1& -1) \\
$16a_{312423}$ &(  1&  1&  0& -1&  0) \\
$16a_{324978}$ &(  1& -2&  0& -1&  0) \\
$16a_{329518}$ &(  0& -2& -1&  1&  0) \\
$16n_{404143}$ &(  0& -1& -1&  1&  1) \\
$16n_{526315}$ &(  1&  0&  1& -1&  0) \\
$16n_{918713}$ &(  0&  0&  1&  2& -1) \\
\end{tabular}
}

\vspace{1cm}
\parbox[t]{5cm}{
\centering
\mbox{} \\
\includegraphics[scale=\scaling]{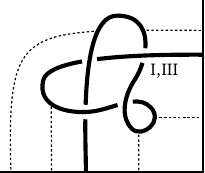} \\
template $t_{26}$ \\
\mbox{} \\
\begin{tabular}{l@{\;=\;}r@{, }r@{, }r@{ $\mid$ }r@{, }r}
$16a_{326826}$  &( 0&  2& -1&  0&  1) \\
$16a_{330763}$  &( 0&  1&  1&  1& -1) \\
$16a_{339566}$  &( 0&  0& -1&  1&  1) \\
$16a_{357272}$  &(-1&  1&  1&  1&  0) \\
\end{tabular}
}
\quad
\parbox[t]{5cm}{
\centering
\mbox{} \\
\includegraphics[scale=\scaling]{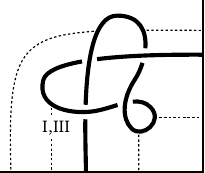} \\
template $t_{27}$ \\
\mbox{} \\
\begin{tabular}{l@{\;=\;}r@{, }r@{, }r@{ $\mid$ }r@{, }r}
$16a_{106900}$  &(  1&  0&  2&  1&  0) \\
$16a_{155718}$  &(  1&  0& -1&  1&  0) \\
$16a_{359219}$  &(  0& -2& -1&  0& -1) \\
\end{tabular}
}

\vspace{1cm}
\parbox[t]{5cm}{
\centering
\mbox{} \\
\includegraphics[scale=\scaling]{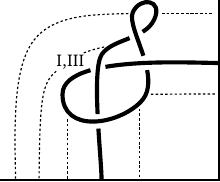} \\
template $t_{28}$ \\
\mbox{} \\
\begin{tabular}{l@{\;=\;}r@{, }r@{, }r@{, }r@{ $\mid$ }r@{, }r}
$14a_{19517}$  &( 1&  0&  0& -1&  1&-1) \\
$16a_{309401}$ &( 2&  0&  0& -1&  1& 0) \\
$16a_{354193}$ &( 1&  0&  0& -2&  1& 0) \\
$16n_{925482}$ &( 2& -1&  0&  2&  1& 0) \\
\end{tabular}
}

\normalsize

\clearpage
\newpage

\vspace{1cm}
\noindent
Christoph Lamm \\ \noindent
R\"{u}ckertstr. 3, 65187 Wiesbaden \\ \noindent
Germany \\ \noindent
e-mail: christoph.lamm@web.de

\end{document}